\documentclass[aos,MSNbibl,dvips]{arximspdf}
\usepackage{mathbh}
\usepackage{graphicx}
\usepackage{breakurl}

\doi{10.1214/13-AOS1142} 
\volume{41}
\issue{5}
\pubyear{2013}
\firstpage{2292}
\lastpage{2323}

\makeatletter
\newcommand{\rrVert}{\Vert}
\newcommand{\rrvert}{\vert}
\newcommand{\llVert}{\Vert}
\newcommand{\llvert}{\vert}
\def\mc#1{\mathcal{#1}}
\def\what#1{\widehat{#1}}

\def\R{\mathbb{R}}
\def\N{\mathbb{N}}

\def\P{\mathbb{P}} 

\newcommand{\argmin}{\operatorname{argmin}}
\newcommand{\diag}{\operatorname{diag}}
\newcommand{\sign}{\operatorname{sign}}

\newcommand{\minimize}{\mathrm{minimize}}

\newtheorem{theorem}{Theorem}
\newtheorem{lemma}{Lemma}
\newtheorem{corollary}{Corollary}
\newtheorem{proposition}{Proposition}
\newproclaim{definition}{Definition}
\newproclaim{assumption}{Assumption}
\newproclaim{condition}{Condition}

\newcommand{\eqref}[1]{(\ref{#1})}
\renewcommand{\citep}[1]{(\citeauthor{#1} \citeyear{#1})}
\def\mc{\mathcal}
\def\surrloss{\varphi} 
\def\dag{G} 
\def\limitgraph{G_\limitlaw} 
\def\preflabel{Y} 
\def\preflabelspace{\mc{Y}}
\def\queryspace{\mc{Q}}
\def\structure{s} 
\def\structurespace{\mc{S}} 
\def\measures{\mc{M}} 
\def\limitlaw{\mu} 
\def\condloss{\ell} 
\def\condsurrloss{\ell_\surrloss} 
\def\E{\mathbb{E}} 
\def\defeq{\triangleq} 
\def\cas{\stackrel{\mathrm{a.s.}}{\rightarrow}} 
\def\defeq{:=}
\def\funclass{\mc{F}} 
\def\risk{R} 
\def\emprisk{\what{R}} 
\newcommand{\NP}{\mathit{NP}}

\def\order{\mc{O}} 
\def\queryexp{\beta}
\def\queryprob{p}
\def\numresults{{m}} 
\def\batch{\mc{B}}
\def\pinv{\dagger}

\def\surrbound{B}
\def\lipbound{L}

\def\cas{\stackrel{\mathrm{a.s.}}{\rightarrow}}
\def\cp{\stackrel{p}{\rightarrow}}
\def\cd{\stackrel{d}{\rightarrow}}

\newcommand{\onevec}{\mathbh{1}}

\makeatother

\begin{document}
\begin{frontmatter}

\title{The asymptotics of ranking algorithms}
\runtitle{The asymptotics of ranking algorithms}

\begin{aug}
\author[a]{\fnms{John C.} \snm{Duchi}\corref{}\thanksref{t1,t2}\ead[label=e1]{jduchi@cs.berkeley.edu}},
\author[b]{\fnms{Lester} \snm{Mackey}\thanksref{t1,t2}\ead[label=e2]{lmackey@stanford.edu}}
\and
\author[a]{\fnms{Michael I.} \snm{Jordan}\thanksref{t2}\ead[label=e3]{jordan@stat.berkeley.edu}}
\runauthor{J. C. Duchi, L. Mackey and M. I. Jordan}
\affiliation{University of California, Berkeley, Stanford University
and University~of~California,~Berkeley}
\address[a]{J. C. Duchi\\
M. I. Jordan\\
Departments of EECS and Statistics\\
University of California, Berkeley\\
Berkeley, California\\
USA\\
\printead{e1} \\
\phantom{E-mail:\ }\printead*{e3}}

\address[b]{L. Mackey\\
Department of Statistics\\
Stanford University\\
Stanford, California\\
USA\\
\printead{e2}}
\thankstext{t1}{Supported by DARPA through the National Defense Science
and Engineering Graduate Fellowship Program (NDSEG).}
\thankstext{t2}{Supported in part by the U.S. Army Research Laboratory and
the U.S. Army Research Office under contract/Grant W911NF-11-1-0391.}
\end{aug}

\received{\smonth{4} \syear{2012}}
\revised{\smonth{1} \syear{2013}}

%
\begin{abstract}
We consider the predictive problem of supervised ranking, where the
task is
to rank sets of candidate items returned in response to queries. Although
there exist statistical procedures that come with guarantees of consistency
in this setting, these procedures require that individuals provide a
complete ranking of all items, which is rarely feasible in practice.
Instead, individuals routinely provide partial preference information, such
as pairwise comparisons of items, and more practical approaches to ranking
have aimed at modeling this partial preference data directly. As we show,
however, such an approach raises serious theoretical challenges.
Indeed, we
demonstrate that many commonly used surrogate losses for pairwise comparison
data do not yield consistency; surprisingly, we show inconsistency even in
low-noise settings. With these negative results as motivation, we
present a
new approach to supervised ranking based on aggregation of partial
preferences, and we develop $U$-statistic-based empirical risk minimization
procedures. We present an asymptotic analysis of these new procedures,
showing that they yield consistency results that parallel those available
for classification. We complement our theoretical results with an
experiment studying the new procedures in a large-scale web-ranking task.
\end{abstract}

%

%
%
%

\end{frontmatter}

\section{Introduction}

Recent years have seen significant developments in the theory of
classification, most notably binary classification, where strong theoretical
results are available that quantify rates of convergence and shed light on
qualitative aspects of the problem \cite{Zhang04b,BartlettJoMc06}.
Extensions to multi-class classification have also been explored, and
connections to the theory of regression are increasingly well
understood, so
that overall a satisfactory theory of supervised machine learning has
begun to
emerge \cite{Zhang04a,Steinwart07}.

In many real-world problems in which labels or responses are available,
however, the problem is not merely to classify or predict a real-valued
response, but rather to list a set of items in order. The theory of
supervised learning cannot be considered complete until it also
provides a
treatment of such \emph{ranking} problems. For example, in information
retrieval, the goal is to rank a set of documents in order of
relevance to a user's search query; in medicine, the object is often to
rank drugs in order of probable curative outcomes for a given disease; and
in recommendation or advertising systems, the aim is to present a set of
products in order of a customer's willingness to purchase or consume.
In each
example, the intention is to order a set of items in accordance with the
preferences of an individual or population. While such problems are often
converted to classification problems for simplicity (e.g., a
document is
classified as ``relevant'' or not), decision makers frequently require the
ranks (e.g., a search engine must display documents in a particular
order on the page). Despite its ubiquity, our statistical understanding of
ranking falls short of our understanding of classification and regression.
Our aim here is to characterize the statistical behavior of computationally
tractable inference procedures for ranking under natural data-generating
mechanisms.

We consider a general decision-theoretic formulation of the \emph{supervised
ranking problem} in which preference data are drawn i.i.d. from an unknown
distribution, where each datum consists of a \emph{query}, $Q \in
\queryspace$,
and a \emph{preference judgment}, $\preflabel\in\preflabelspace$,
over a set
$m$ of candidate items that are available based on the query $Q$. The exact
nature of the query and preference judgment depend on the ranking
context. In the setting of information retrieval, for example, each datum
corresponds to a user issuing a natural language query and expressing a
preference by selecting or clicking on zero or more of the returned results.
The statistical task is to discover a function that provides a
query-specific ordering of items that best respects the observed preferences.
This query-indexed setting is especially natural for tasks like information
retrieval in which a different ranking of webpages is needed for each natural
language query.

Following existing literature, we estimate a \emph{scoring function}
$f\dvtx \queryspace\rightarrow\R^{\numresults}$, where $f(q)$ assigns a
score to
each of $\numresults$ candidate items for the query $q$, and the
results are
ranked according to their scores \cite{HerbrichGrOb00,FreundIyScSi03}.
Throughout the paper, we adopt a decision-theoretic perspective and assume
that given a query-judgment pair $(Q, \preflabel)$, we evaluate the scoring
function $f$ via a loss $L(f(Q), \preflabel)$. The goal is to choose the
$f$ minimizing the risk
%
%
\begin{equation}
\label{eqnrisk} \risk(f) \defeq\E\bigl[L\bigl(f(Q), \preflabel\bigr)\bigr].
\end{equation}
While minimizing the risk \eqref{eqnrisk} directly is in general intractable,
researchers in machine learning and information retrieval have developed
surrogate loss functions that yield procedures for selecting $f$. Unfortunately,
as we show, extant procedures fail to solve the ranking problem under reasonable
data generating mechanisms. The goal in the remainder of the paper is
to explain
this failure and to propose a novel solution strategy based on preference
aggregation.

Let us begin to elucidate the shortcomings of current approaches to ranking.
One main problem lies in their unrealistic assumptions about available data.
The losses proposed and most commonly used for evaluation in the information
retrieval literature \cite{ManningRaSc08,JarvelinKe04} have a common form,
generally referred to as (Normalized) Discounted Cumulative Gain ((N)DCG).
The NDCG family requires that the preference judgments $\preflabel$
associated with the datum $(Q, \preflabel)$ be a vector
$\preflabel\in\R^\numresults$ of \emph{relevance scores} for the
entire set of items; that is, $\preflabel_j$ denotes the real-valued relevance
of item $j$ to the query $Q$. While having complete preference information
makes it possible to design procedures that asymptotically minimize
NDCG losses (e.g., \cite{CossockZh08}), in practice such complete
preferences are unrealistic: they are expensive to collect and
difficult to trust. In biological applications, evaluating the effects
of all drugs involved in a study---or all doses---on a single subject
is infeasible. In web search, users click on only one or two results: no
feedback is available for most items. Even when practical and ethical
considerations do not preclude collecting complete preference information
from participants in a study, a long line of psychological work
highlights the inconsistency with which humans assign numerical values
to multiple objects (e.g., \cite{ShiffrinNo94,StewartBrCh05,Miller56}).

The inherent practical difficulties that arise in using losses based on
relevance scores has led other researchers to propose loss functions
that are suitable for \emph{partial preference data} \cite{Joachims02,FreundIyScSi03,DekelMaSi03}. Such data arise naturally in a number of
real-world situations; for example, a patient's prognosis may improve
or deteriorate after administration of treatment, competitions and
sporting matches provide paired results, and shoppers at a store purchase
one item but not others. Moreover, the psychological literature shows
that human beings are quite good at performing pairwise distinctions and
forming relative judgments (see, e.g., \cite{Saaty08} and references therein).

More formally, let $\alpha\defeq f(Q) \in\R^\numresults$ denote the vector
of predicted scores for each item associated with query $Q$. If a preference
$\preflabel$ indicates that item $i$ is preferred to $j$ then the natural
associated loss is the zero-one loss $L(\alpha, \preflabel) = 1
({\alpha_i \le\alpha_j} )$. Minimizing such a loss is well
known to be computationally
intractable; nonetheless, the classification
literature \cite{Zhang04a,Zhang04b,BartlettJoMc06,Steinwart07} has shown
that it is possible to design convex Fisher-consistent surrogate losses for
the 0--1 loss in classification settings and has linked Fisher
consistency to consistency. By reduction to classification, similar
consistency results are possible in certain bipartite or binary ranking
scenarios \cite{ClemenconLuVa08}. One might therefore hope to make use of
these surrogate losses in the ranking setting to obtain similar guarantees.
Unfortunately, however, this hope is not borne out; as we illustrate in
Section~\ref{secconsistency}, it is generally computationally
intractable to
minimize any Fisher-consistent loss for ranking, and even in favorable
low-noise cases, convex surrogates that yield Fisher consistency for binary
classification fail to be Fisher-consistent for ranking.

We find ourselves at an impasse: existing methods based on practical
data-collection strategies do not yield a satisfactory theory, and
those methods that do have theoretical justification are not practical.
Our approach to this difficulty is to take a new approach to supervised
ranking problems in which partial preference data are aggregated before
being used for estimation. The point of departure for this approach
is the notion of \emph{rank aggregation} (e.g., \cite{DworkKuNaSi01}),
which has a long history in voting \cite{Borda1781}, social choice
theory \cite{Condorcet1785,Arrow51} and statistics \cite{Thurstone27,Mallows57}.
In Section~\ref{secaggregation}, we discuss some of the
ways in which partial preference data can be aggregated, and we propose a
new family of $U$-statistic-based loss functions that are computationally
tractable. Sections~\ref{secconsistency} and \ref{secu-statistics}
present a theoretical analysis of procedures based on these loss functions,
establishing their consistency. We provide a further discussion of
practical rank aggregation strategies in Section~\ref{secrank-aggregation}
and present experimental results in Section~\ref{secexperiments}.
Section~\ref{secconclusions} contains
our conclusions, with proofs deferred to appendices.

\section{Ranking with rank aggregation}
\label{secaggregation}

We begin by considering several ways in which partial preference data
arise in practice. We then turn to a formal treatment of our aggregation-based
strategy for supervised ranking.
\begin{longlist}[1.]
\item[1.]
\label{itempaired-comparison}
\emph{Paired comparison data}. Data in which an individual judges one item
to be preferred over another in the context of a query are common.
Competitions and sporting matches, where each pairwise comparison may be
accompanied by a magnitude such as a difference of scores,
naturally generate such data. In practice, a single individual will not
provide feedback for all possible pairwise comparisons, and we do not assume
transitivity among the observed preferences for an individual. Thus, it is
natural to model the pairwise preference judgment space
$\preflabelspace$ as
the set of weighted directed graphs on $\numresults$ nodes.
\item[2.]\emph{Selection data}. A ubiquitous source of partial
preference information is the selection behavior of a user presented
with a
small set of potentially ordered items. For example, in response to a
search query, a web search engine presents an ordered list of webpages and
records the URL a user clicks on, and a store records inventory and tracks
the items customers purchase. Such selections provide partial
information: that a user or customer prefers one item to others presented.
\item[3.]\emph{Partial orders}. An individual may also provide preference
feedback in terms of a partial ordering over a set of candidates or items.
In the context of elections, for example,
each preference judgment $\preflabel\in\preflabelspace$ specifies a partial
order $\prec_\preflabel$ over candidates such that candidate $i$ is
preferred to candidate $j$ whenever $i \prec_\preflabel j$.
A~partial order need not specify a preference between every pair of items.
\end{longlist}

Using these examples as motivation, we wish to develop a formal
treatment of
ranking based on aggregation. To provide intuition for the framework
presented in the remainder of this section, let us consider a simple
aggregation strategy appropriate for the case of paired comparison
data. Let
each relevance judgment $\preflabel\in\preflabelspace$ be a weighted
adjacency matrix where the $(i, j)$th entry expresses a preference for item
$i$ over $j$ whenever this entry is nonzero. In this case, a natural
aggregation strategy is to average all observed adjacency matrices for
a fixed
query. Specifically, for a set of adjacency matrices
$\{\preflabel_l\}_{l=1}^k$ representing user preferences for a given
query, we
form the average $(1/k) \sum_{l=1}^k \preflabel_l$. As $k \rightarrow
\infty$,
the average adjacency matrix captures the mean population preferences,
and we
thereby obtain complete preference information over the $\numresults$ items.

This averaging of partial preferences is one example of a general class of
aggregation strategies that form the basis of our theoretical framework.
To formalize this notion, we modify the loss formulation slightly and
hereafter assume that the loss function $L$ is a mapping $\R
^\numresults
\times\structurespace\rightarrow\R$, where $\structurespace$ is a
problem-specific \emph{structure space}. We further assume the
existence of a
series of \emph{structure functions}, $\structure_k \dvtx \preflabelspace^k
\rightarrow\structurespace$, that map sets of preference judgments
$\{\preflabel_j\}$ into $\structurespace$.
The loss $L$ depends on the preference feedback
$(\preflabel_1,\ldots,\preflabel_k)$ for a given query only via the structure
$s_k(\preflabel_1,\ldots,\preflabel_k)$. In the example of the previous
paragraph, $\structurespace$ is the set of $\numresults\times
\numresults$
adjacency matrices, and $\structure_k(\preflabel_1,\ldots,
\preflabel_k) =
(1/k) \sum_{l=1}^k \preflabel_l$. A typical loss for this setting is the
pairwise loss \cite{FreundIyScSi03,Joachims02}
\[
L\bigl(\alpha, \structure(\preflabel_1, \ldots,
\preflabel_k)\bigr) \equiv L(\alpha, A) \defeq\sum
_{i < j} A_{ij} 1 ({\alpha_i \le
\alpha_j} ) + \sum_{i > j} A_{ij}
1 ({\alpha_i < \alpha_j} ),
\]
where $\alpha$ is a set of scores and $A =
\structure_k(\preflabel_1, \ldots, \preflabel_k)$ is the average adjacency
matrix with entries $A_{ij}$. In Section~\ref{secrank-aggregation},
we provide other examples of structure functions for different
data collection mechanisms and losses. Hereafter, we abbreviate
$\structure_k(\preflabel_1, \ldots, \preflabel_k)$ as
$\structure(\preflabel_1, \ldots, \preflabel_k)$ whenever the input
length $k$
is clear from context.

To meaningfully characterize the asymptotics of inference procedures, we
make a mild assumption on the limiting behavior of the structure functions.
%
\renewcommand{\theassumption}{\Alph{assumption}}
\begin{assumption}
\label{assumptionlimit-structures}
Fix a query $Q = q$. Let the sequence $\preflabel_1,\preflabel_2,
\ldots$
be drawn i.i.d. conditional on $q$, and define the random
variables $S_k \defeq\structure(\preflabel_1, \ldots,
\preflabel_k)$. If $\limitlaw^k_q$ denotes the distribution of $S_k$,
there exists a limiting law $\limitlaw_q$ such that
\[
\limitlaw^k_q \cd\limitlaw_q \qquad\mbox{as } k
\rightarrow\infty.
\]
\end{assumption}

For example, the averaging structure function satisfies
Assumption~\ref{assumptionlimit-structures} so long as $\E
[|\preflabel_{ij}|
\mid Q] < \infty$ with probability 1.
Aside from the requirements of
Assumption~\ref{assumptionlimit-structures}, we allow arbitrary
aggregation within the structure function.

In addition, our main assumption on the loss function $L$ is as follows:
%
\begin{assumption}
\label{assumptionloss-continuity}
The loss function $L \dvtx \R^\numresults\times\structurespace
\rightarrow\R$
is bounded in $[0, 1]$, and, for any fixed vector $\alpha\in
\R^\numresults$, $L(\alpha, \cdot)$ is continuous in the topology of
$\structurespace$.
\end{assumption}

With our assumptions on the asymptotics of the structure function
$\structure$ and the loss $L$ in place, we now describe the risk
functions that guide our design of inference procedures. We begin with
the pointwise conditional risk, which maps predicted scores and a measure
$\limitlaw$ on $\structurespace$ to $[0, 1]$:
%
%
\begin{equation}
\condloss\dvtx \R^m \times\measures(\structurespace) \rightarrow[0,
1] \qquad\mbox{where } \condloss(\alpha, \limitlaw) \defeq \int L(\alpha, \structure) \,d
\limitlaw(\structure). \label{eqncondloss}
\end{equation}
Here, $\measures(\structurespace)$ denotes the closure of the subset of
probability measures on the set $\structurespace$ for which $\condloss
$ is
defined.
For any query $q$ and $\alpha\in\R^m$, we have $\lim_i \condloss
(\alpha,
\limitlaw_q^i) = \condloss(\alpha, \limitlaw_q)$ by the definition of
convergence in distribution. This convergence motivates our decision-theoretic
approach.

Our goal in ranking is thus to minimize the risk
%
%
\begin{equation}
\label{eqnrisk-def} \risk(f) \defeq \sum_q
\queryprob_q \condloss\bigl(f(q), \limitlaw_q\bigr),
\end{equation}
where $\queryprob_q$ denotes the probability that the query $Q = q$ is issued.
The risk of the scoring function $f$ can also be obtained in the limit
as the
number of preference judgments for each query goes to infinity:
%
%
\begin{equation}
\label{eqnalternate-risk} \risk(f) = \lim_k \E \bigl[ L
\bigl(f(Q), \structure(\preflabel_1, \ldots, \preflabel_k)
\bigr) \bigr] = \lim_k \sum_q
\queryprob_q \condloss\bigl(f(q), \limitlaw_q^k
\bigr).
\end{equation}
That the limiting expectation \eqref{eqnalternate-risk} is equal
to the risk \eqref{eqnrisk-def} follows from the definition of
weak convergence.

We face two main difficulties in the study of the minimization of the
risk \eqref{eqnrisk-def}. The first difficulty is that of \emph{Fisher
consistency} mentioned previously: since $L$ may be nonsmooth in the
function $f$ and is typically intractable to minimize, when will the
minimization of a tractable surrogate lead to the minimization of the
loss \eqref{eqnrisk-def}? We provide a precise formulation of and
answer to
this question in Section~\ref{secconsistency}.\vadjust{\goodbreak} In addition, we
demonstrate the inconsistency of many commonly used pairwise ranking
surrogates and show that aggregation leads to tractable Fisher consistent
inference procedures for both complete and partial data losses.

The second difficulty is that of \emph{consistency}: for a given Fisher
consistent surrogate for the risk \eqref{eqnrisk-def}, are there tractable
statistical procedures that converge to a minimizer of the risk? Yes:
in Section~\ref{secu-statistics}, we
develop a new family of aggregation losses based on $U$-statistics of
increasing order, showing that uniform laws of large numbers hold for the
resulting $M$-estimators.

\section{Fisher consistency of surrogate risk minimization}
\label{secconsistency}

In this section, we formally define the Fisher consistency of a
surrogate loss
and give general necessary and sufficient conditions for consistency to hold
for losses satisfying Assumption~\ref{assumptionloss-continuity}. To begin,
we assume that the space $\queryspace$ of queries is countable (or
finite) and thus
bijective with $\N$. Recalling the definition \eqref{eqnrisk-def} of
the risk
and the pointwise conditional risk \eqref{eqncondloss}, we define the
Bayes risk for
$\risk$ as the minimal risk over all measurable functions $f \dvtx
\queryspace
\rightarrow\R^\numresults$:
\[
\risk^* \defeq\inf_f \risk(f) = \sum
_q \queryprob_q \inf_{\alpha\in\R^\numresults}
\condloss(\alpha, \limitlaw_q).
\]
The second equality follows because $\queryspace$ is countable and
the infimum is taken over all measurable functions.

Since it is infeasible to minimize the risk \eqref{eqnrisk-def}
directly, we
consider a bounded-below surrogate $\surrloss$ to minimize in place of
$L$. For each structure $\structure\in\structurespace$, we write
$\surrloss(\cdot, \structure) \dvtx \R^\numresults\rightarrow\R_+$,
and we
assume that for $\alpha\in\R^\numresults$, the function $\structure
\mapsto
\surrloss(\alpha, \structure)$ is continuous with respect to the
topology on
$\structurespace$. We then define the conditional $\surrloss$-risk as
%
%
\begin{equation}
\label{eqncond-surrogate} \condsurrloss(\alpha, \limitlaw) \defeq \int
_{\structurespace} \surrloss(\alpha, \structure) \,d\limitlaw (\structure)
\end{equation}
and the asymptotic $\surrloss$-risk of the function $f$ as
%
%
\begin{equation}
\label{eqnsurrogate-risk} \risk_\surrloss(f) \defeq\sum
_q \queryprob_q \condsurrloss\bigl(f(q),
\limitlaw_q\bigr),
\end{equation}
whenever each $\condsurrloss(f(q), \limitlaw_q)$ exists [otherwise
$\risk_\surrloss(f) = +\infty$]. The optimal $\surrloss$-risk is
defined to be
$\risk^*_\surrloss\defeq\inf_f \risk_\surrloss(f)$, and
throughout we make
the assumption that there exist measurable $f$ such that $\risk
_\surrloss(f) <
+\infty$ so that $\risk_\surrloss^*$ is finite. The following is our general
notion of Fisher consistency.
%
\begin{definition}
\label{definitionconsistency}
The surrogate loss $\surrloss$ is \emph{Fisher-consistent} for the
loss $L$
if for any $\{\queryprob_q\}$ and probability measures $\limitlaw_q
\in
\measures(\structurespace)$, the convergence
\[
\risk_\surrloss(f_n) \rightarrow\risk_\surrloss^*\qquad
\mbox{implies } \risk(f_n) \rightarrow\risk^*.
\]
\end{definition}
To achieve more actionable risk bounds and to more accurately compare surrogate
risks, we also draw upon a uniform statement of consistency:\vadjust{\goodbreak}
%
\begin{definition}
\label{definitionuniform-consistency}
The surrogate loss $\surrloss$ is \emph{uniformly Fisher-consistent}
for the
loss $L$ if for any $\varepsilon> 0$, there exists a $\delta(\varepsilon)
> 0$
such that for any $\{\queryprob_q\}$ and probability measures
$\limitlaw_q
\in\measures(\structurespace)$,
%
%
\begin{equation}
\risk_\surrloss(f) < \risk_\surrloss^* + \delta(\varepsilon)
\qquad\mbox{implies } \risk(f) < \risk^* + \varepsilon. \label{eqnuniform-consistency}
\end{equation}
\end{definition}

The bound \eqref{eqnuniform-consistency} is equivalent to the
assertion that there exists a nondecreasing function $\zeta$ such
that $\zeta(0) = 0$ and $\risk(f) - \risk^* \le\zeta(\risk
_\surrloss(f) -
\risk_\surrloss^*)$. Bounds of this form have been completely
characterized in
the case of binary classification~\cite{BartlettJoMc06},
and Steinwart~\cite{Steinwart07} has given necessary and sufficient conditions for
uniform Fisher-consistency to hold in general
risk minimization problems. We now turn to analyzing conditions under
which a surrogate loss $\surrloss$ is Fisher-consistent for ranking.

\subsection{General theory}
\label{secgeneral-consistency}

The main approach in establishing conditions for the surrogate risk
Fisher consistency in Definition~\ref{definitionconsistency} is to
move from global
conditions for Fisher consistency to local, pointwise Fisher consistency.
Following the treatment of Steinwart~\cite{Steinwart07}, we begin by defining a
function measuring the discriminating ability of the surrogate
$\surrloss$:
%
%
\begin{equation}
H(\varepsilon) \defeq\inf_{\limitlaw\in\measures(\structurespace),
\alpha} \Bigl\{\condsurrloss(\alpha,
\limitlaw) - \inf_{\alpha'}\condsurrloss\bigl(\alpha',
\limitlaw\bigr) \mid \condloss(\alpha, \limitlaw) - \inf_{\alpha'}
\condloss\bigl(\alpha', \limitlaw\bigr) \ge\varepsilon \Bigr\}.
\label{eqnsuboptimality-def}
\end{equation}
This function is familiar from work on surrogate risk Fisher
consistency in
classification \cite{BartlettJoMc06} and measures surrogate risk
suboptimality as a function of risk suboptimality.
A reasonable conditional $\surrloss$-risk will declare a set of scores
$\alpha\in\R^\numresults$ suboptimal whenever the conditional risk
$\condloss$ declares them suboptimal.
This corresponds to $H(\varepsilon) > 0$ whenever $\varepsilon> 0$, and
we call any loss satisfying this condition \emph{pointwise consistent}.

From these definitions, we can conclude the following consistency result,
which is analogous to the results of \cite{Steinwart07}. For
completeness, we
provide a proof in
the supplementary material \cite{DuchiMaJo13supp}.
%
\begin{proposition}
\label{propositionuniform-consistency}
Let $\surrloss\dvtx  \R^\numresults\times\structurespace\rightarrow\R
_+$ be
a bounded-below loss function such that for some $f$, $\risk_\surrloss
(f) <
+\infty$. Then $\surrloss$ is pointwise consistent if and only if the
uniform Fisher-consistency definition \eqref{eqnuniform-consistency} holds.
\end{proposition}
%

\newcommand{\partition}{\mc{A}} 

Proposition~\ref{propositionuniform-consistency} makes it clear that
pointwise consistency for general measures $\limitlaw$ on the set of structures
$\structurespace$ is a stronger condition than that of Fisher
consistency in
Definition~\ref{definitionconsistency}. In some situations, however,
it is
possible to connect the weaker surrogate risk Fisher consistency of
Definition~\ref{definitionconsistency} with uniform Fisher
consistency and pointwise
consistency. Ranking problems with appropriate choices of the space
$\structurespace$ give rise to such connections. Indeed, consider the
following:
%
\begin{assumption}
\label{assumptionfiniteness}
The space of possible structures $\structurespace$ is finite, and the loss
$L$ is discrete, meaning that it takes on only finitely many values.
\end{assumption}

Binary and multiclass classification provide examples of settings in which
Assumption~\ref{assumptionfiniteness} is appropriate, since the set
of structures
$\structurespace$ is the set of class labels, and $L$ is usually a
version of
the $0$--$1$ loss. We also sometimes make a weaker version of
Assumption~\ref{assumptionfiniteness}:
%
\renewcommand{\theassumption}{\Alph{assumption}$^{\prime}$}
\setcounter{assumption}{2}
\begin{assumption}
\label{assumptioncompact-S}
The (topological) space of possible structures $\structurespace$ is compact,
and for some $d \in\N$
there exists a partition $\partition_1, \ldots, \partition_d$ of
$\R^\numresults$ such that for any $\structure\in\structurespace$,
\[
L(\alpha, \structure) = L\bigl(\alpha', \structure\bigr)\qquad
\mbox{whenever } \alpha, \alpha' \in\partition_i.
\]
\end{assumption}

Assumption~\ref{assumptioncompact-S} may be more natural in ranking settings
than Assumption~\ref{assumptionfiniteness}. The
compactness assumption holds, for example, if $\structurespace
\subset\R^\numresults$ is closed and bounded,
such as in our pairwise aggregation example in Section~\ref{secaggregation}.
Losses $L$ that depend only on the relative order of the coordinate
values of $\alpha\in\R^\numresults$---common in ranking problems---provide
a collection of examples for which the partitioning condition holds.


Under Assumption~\ref{assumptionfiniteness} or \ref
{assumptioncompact-S}, we
can provide a definition of local consistency that is often more user-friendly
than pointwise consistency \eqref{eqnsuboptimality-def}:
%
\begin{definition}
\label{definitionstructure-consistent}
Let $\surrloss$ be a bounded-below surrogate loss such that
$\surrloss(\cdot, \structure)$ is continuous for all $\structure
\in\structurespace$. The function $\surrloss$ is
\emph{structure-consistent} with respect to the loss $L$ if for all
$\limitlaw\in\measures(\structurespace)$,
\[
\label{eqnstructure-consistent} \condsurrloss^*(\limitlaw) \defeq \inf
_\alpha\condsurrloss(\alpha, \limitlaw) < \inf
_\alpha \Bigl\{ \condsurrloss(\alpha, \limitlaw) \mid\alpha
\notin\mathop{\argmin}_{\alpha'} \condloss\bigl(\alpha', \limitlaw\bigr)
\Bigr\}.
\]
\end{definition}

Definition~\ref{definitionstructure-consistent} describes the set
of loss functions $\surrloss$ satisfying the intuitively desirable property
that the surrogate $\surrloss$ cannot be minimized if the scores
$\alpha\in\R^\numresults$ are restricted to not minimize the loss $L$.
As we see presently, Definition~\ref{definitionstructure-consistent} captures
exactly what it means for a surrogate loss $\surrloss$ to be Fisher-consistent
when one of Assumptions \ref{assumptionfiniteness}
or \ref{assumptioncompact-S} holds. Moreover, the set of Fisher-consistent
surrogates coincides with the set of uniformly Fisher-consistent
surrogates in
this case. The following theorem formally states this result; we give a proof
in
the supplementary material \cite{DuchiMaJo13supp}.
%
\begin{theorem}
\label{theoremstructure-consistency}
Let $\surrloss\dvtx  \R^m \times\structurespace\rightarrow\R_+$ satisfy
$\risk_\surrloss(f) < +\infty$ for some measurable $f$.
If Assumption~\ref{assumptionfiniteness} holds, then:
\begin{longlist}[(a)]
\item[(a)]\label{theoremstructure-consistency-a}
If $\surrloss$ is structure consistent
(Definition~\ref{definitionstructure-consistent}), then $\surrloss$ is
uniformly Fisher-consistent for the loss $L$
(Definition~\ref{definitionuniform-consistency}).
\item[(b)]\label{theoremstructure-consistency-b}
If $\surrloss$ is Fisher-consistent for the loss $L$
(Definition~\ref{definitionconsistency}), then $\surrloss$ is structure
consistent.
\end{longlist}
If the function
$\surrloss(\cdot, \structure)$ is convex for $\structure
\in\structurespace$, and for $\limitlaw\in\measures
(\structurespace)$
the conditional risk $\condsurrloss(\alpha, \limitlaw) \rightarrow
\infty$
as $\llVert {\alpha}\rrVert  \rightarrow\infty$, then
Assumption~\ref{assumptioncompact-S} implies
\textup{(a)} and
\textup{(b)}.
\end{theorem}

Theorem~\ref{theoremstructure-consistency} shows that as long as
Assumption~\ref{assumptionfiniteness} holds, pointwise consistency, structure
consistency, and both uniform and nonuniform surrogate loss consistency
coincide. These four also coincide under the weaker
Assumption~\ref{assumptioncompact-S} so long as the surrogate is
$0$-coercive, which is not restrictive in practice.
As a final note, we recall a result due to Steinwart \cite{Steinwart07},
which gives general necessary and sufficient conditions for the
consistency in
Definition~\ref{definitionconsistency} to hold, using a weaker
version of the suboptimality function \eqref{eqnsuboptimality-def} that
depends on $\limitlaw$:
%
%
\begin{equation}
\label{eqnlocal-suboptimality-def} H(\varepsilon, \limitlaw) \defeq\inf
_{\alpha} \Bigl\{ \condsurrloss(\alpha, \limitlaw) - \inf
_{\alpha'} \condsurrloss\bigl(\alpha', \limitlaw\bigr)
\mid\condloss(\alpha, \limitlaw) - \inf_{\alpha'} \condloss \bigl(
\alpha', \limitlaw\bigr) \ge\varepsilon \Bigr\}.
\end{equation}

\begin{proposition}[(Steinwart \cite{Steinwart07}, Theorems 2.8 and 3.3)]
\label{propositionall-consistency}
The suboptimality function \eqref{eqnlocal-suboptimality-def}
satisfies $H(\varepsilon, \limitlaw_q) > 0$ for any
$\varepsilon> 0$ and $\limitlaw_q$ with $q \in\queryspace$ and
$\queryprob_q
> 0$ if and only if $\surrloss$ is
Fisher-consistent for the loss $L$ (Definition~\ref{definitionconsistency}).
\end{proposition}

As a corollary of this result, any structure-consistent
surrogate loss $\surrloss$ (in the sense of
Definition~\ref{definitionstructure-consistent}) is Fisher-consistent
for the loss
$L$ whenever the conditional risk $\condloss(\alpha, \limitlaw)$ has finite
range, so that $\alpha\notin\argmin_{\alpha'} \condloss(\alpha',
\limitlaw) \neq\varnothing$ implies the existence of an $\varepsilon> 0$ such
that $\condloss(\alpha, \limitlaw) - \inf_{\alpha'} \condloss
(\alpha,
\limitlaw) \ge\varepsilon$.

\subsection{The difficulty of Fisher consistency for ranking}
\label{secconsistency-hard}

We now turn to the question of
whether there exist structure-consistent ranking losses. In a preliminary
version of this work \cite{DuchiMaJo10}, we focused on the practical setting
of learning from pairwise preference data and demonstrated that many popular
ranking surrogates are inconsistent for standard pairwise ranking
losses. We
review and generalize our main inconsistency results here, noting that while
the losses considered use pairwise preferences, they perform no
aggregation. Their theoretically poor performance provides motivation
for the
aggregation strategies proposed in this work; we explore the
connections in
Section~\ref{secrank-aggregation} (focusing on pairwise losses in
Section~\ref{secstructured-aggregation}). We provide proofs of our
inconsistency results in
the supplementary material \cite{DuchiMaJo13supp}.

To place ourselves in the general structural setting of the paper, we consider
the structure function $\structure(\preflabel_1, \ldots, \preflabel
_k) =
\preflabel_1$ which performs no aggregation for all~$k$, and we let
$\preflabel$ denote the weighted adjacency matrix of a directed
acyclic graph
(DAG) $\dag$, so that $\preflabel_{ij}$ is the weight of the directed edge
$ ({i} \rightarrow{j} )$ in the graph $\dag$. We consider
a pairwise
loss that imposes a separate penalty for each misordered pair of results:
%
%
\begin{equation}
\label{eqnedge-loss} L(\alpha, \preflabel) = \sum_{i < j}
\preflabel_{ij} 1 ({\alpha _i \le\alpha_j} ) +
\sum_{i > j} \preflabel_{ij} 1 ({
\alpha_i < \alpha_j} ),
\end{equation}
where we distinguish the cases $i < j$ and $i > j$ to avoid doubly
penalizing $1 ({\alpha_i=\alpha_j} )$. When pairwise preference
judgments are
available, use of such losses is common. Indeed, this loss generalizes the
disagreement error described by Dekel et~al. \cite{DekelMaSi03} and is similar to losses
used by Joachims \cite{Joachims02}. If we define $\preflabel_{ij}^\limitlaw
\defeq
\int\preflabel_{ij} \,d\limitlaw(\preflabel)$, then
%
%
\begin{equation}
\condloss(\alpha, \limitlaw) = \sum_{i < j}
\preflabel_{ij}^\limitlaw1 ({\alpha_i \le\alpha
_j} ) + \sum_{i > j} \preflabel_{ij}^\limitlaw1
({\alpha_i < \alpha _j} ). \label{eqnedge-condloss}
\end{equation}
We assume that the number of nodes in any graph $\dag$ (or,
equivalently, the
number of results returned by any query) is bounded by a finite
constant $\numresults$.
Hence, the conditional risk \eqref{eqnedge-condloss} has a finite range;
if there are a finite number of preference labels $\preflabel$ or
the set of weights is compact,
Assumptions \ref{assumptionfiniteness} or \ref{assumptioncompact-S}
are satisfied, whence Theorem~\ref{theoremstructure-consistency} applies.

\subsubsection{General inconsistency}

Let the set $P$ denote the complexity class of problems solvable in polynomial
time and $\NP$ denote the class of nondeterministic polynomial time
problems (see, e.g., \cite{HopcroftUl79}). Our first inconsistency result
(see also~\cite{DuchiMaJo10}, Lemma~7)
is that unless $P = \NP$ (a widely doubted proposition), any loss that is
tractable to minimize cannot be a Fisher-consistent surrogate for the
loss \eqref{eqnedge-loss} and its associated risk.
%
\begin{proposition}
\label{propositionranking-np-hard}
Finding an $\alpha$ minimizing $\condloss$ is $\NP$-hard.
\end{proposition}

In particular, most convex functions are minimizable to an accuracy of
$\varepsilon$ in time polynomial in the dimension of the problem times a multiple
of $\log\frac{1}{\varepsilon}$, known as poly-logarithmic
time \cite{Ben-TalNe01}. Since any $\alpha$ minimizing
$\condsurrloss(\alpha,
\limitlaw)$ must minimize $\condloss(\alpha, \limitlaw)$ for a
Fisher-consistent
surrogate $\surrloss$, and $\condloss(\cdot, \limitlaw)$ has a
finite range
(so that optimizing $\condsurrloss$ to a fixed $\varepsilon$ accuracy is
sufficient), convex surrogate losses are inconsistent for the pairwise
loss \eqref{eqnedge-loss} unless $P = \NP$.

\subsubsection{Low-noise inconsistency}
\label{seclow-noise-inconsistency}

We now turn to showing that, surprisingly, many common convex
surrogates are
inconsistent even in low-noise settings in which it is easy to find an
$\alpha$ minimizing $\condloss(\alpha, \limitlaw)$. (Weaker
versions of the
results in this section appeared in our preliminary paper \cite{DuchiMaJo10}.)
Inspecting the loss definition~\eqref{eqnedge-loss}, a natural choice
for a
surrogate loss is one of the
form \cite{HerbrichGrOb00,FreundIyScSi03,DekelMaSi03}
%
%
\begin{equation}
\surrloss(\alpha, \preflabel) = \sum_{i,j} h(
\preflabel_{ij}) \phi(\alpha_i - \alpha_j),
\label{eqnsurr-edge-loss}
\end{equation}
where $\phi\ge0$ is a convex function, and $h$ is a some function of the
penalties $\preflabel_{ij}$. This surrogate implicitly uses the structure
function $\structure(\preflabel_1, \ldots, \preflabel_k) =
\preflabel_1$ and
performs no preference aggregation. The conditional surrogate risk is thus
$\condsurrloss(\alpha, \limitlaw) = \sum_{i\neq j} h_{ij} \phi
(\alpha_i -
\alpha_j)$, where $h_{ij} \defeq\int h(\preflabel_{ij})
\,d\limitlaw(\preflabel)$. Surrogates of the form \eqref{eqnsurr-edge-loss}
are convenient in margin-based binary classification, where
the complete description by Bartlett, Jordan and McAuliffe \cite{BartlettJoMc06} shows $\phi$ is
Fisher-consistent if and only if it is differentiable at 0 with $\phi'(0)
< 0$.

%

\begin{figure}

\includegraphics{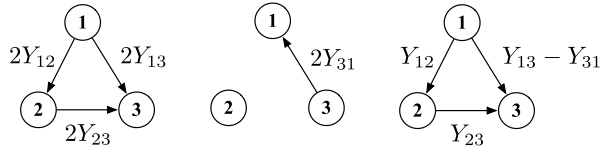}

\caption{The two leftmost DAGs
occur with probability $\frac{1}{2}$, yielding the difference graph
$\limitgraph$ at right, assuming $\preflabel_{23} > \preflabel_{32}$.}
\label{figdag-examples}
\end{figure}

We now precisely define our low-noise setting. For any measure
$\limitlaw$ on
a space~$\preflabelspace$ of adjacency matrices, let the directed graph
$\limitgraph$ be the \emph{difference graph}, that is, the graph with edge
weights $\max\{\preflabel_{ij}^\limitlaw- \preflabel_{ji}^\limitlaw
, 0\}$ on
edges $ ({i} \rightarrow{j} )$, where $\preflabel
_{ij}^\limitlaw= \int
\preflabel_{ij} \,d\limitlaw(\preflabel)$. Then we say that the edge
$ ({i} \rightarrow{j} ) \notin\limitgraph$ if
$\preflabel_{ij}^\limitlaw\le
\preflabel_{ji}^\limitlaw$ (see Figure~\ref{figdag-examples}). We
define the
following low-noise condition based on self-reinforcement of edges in the
difference graph.

%
\begin{definition}
\label{defrev-triangle}
The measure $\limitlaw$ on a set $\preflabel$ of adjacency matrices is
\emph{low-noise} when the corresponding difference graph $\limitgraph$
satisfies the following reverse triangle inequality: whenever there is an
edge $ ({i} \rightarrow{j} )$ and an edge $ ({j}
\rightarrow{k} )$ in $\limitgraph$, the
weight $\preflabel_{ik}^\limitlaw- \preflabel_{ki}^\limitlaw$ on
the edge
$ ({i} \rightarrow{k} )$ is greater than or equal to the
path weight
$\preflabel^\limitlaw_{ij} - \preflabel^\limitlaw_{ji} + \preflabel
^\limitlaw_{jk} -
\preflabel^\limitlaw_{kj}$ on the path $(i\to j\to k)$.
\end{definition}
If $\limitlaw$ satisfies Definition~\ref{defrev-triangle}, its
difference graph $\limitgraph$ is a DAG. Indeed, the definition
ensures that
all global preference information in $\limitgraph$ (the sum of weights along
any path) conforms with and reinforces local preference information (the
weight on a single edge). Hence, we would expect any reasonable ranking
method to be consistent in this setting. Nevertheless, typical pairwise
surrogate losses are inconsistent in this low-noise setting (see also the
weaker Theorem~11 in our preliminary work \cite{DuchiMaJo10}):
%
\begin{theorem}
\label{theoreminconsistency-convex}
Let $\surrloss$ be a loss of the form \eqref{eqnsurr-edge-loss}
and assume $h(0) = 0$. If $\phi$ is convex,
then even in the low-noise
setting of Definition~\ref{defrev-triangle} the loss
$\surrloss$ is not structure-consistent.
\end{theorem}


Given the difficulties we encounter using losses of the
form \eqref{eqnsurr-edge-loss}, it is reasonable to consider a reformulation
of the surrogate. A natural alternative is a margin-based loss, which encodes
a desire to separate ranking scores by large margins dependent on the
preferences in a graph. Similar losses have been proposed, for example,
by~\cite{ShashuaLe02}. The next result shows that convex margin-based
losses are also inconsistent, even in low-noise settings. (See also the weaker
Theorem~12 of our preliminary work \cite{DuchiMaJo10}.)

\begin{theorem}
\label{theoreminconsistency-margin}
Let $h \dvtx \R\rightarrow\R$ and
$\surrloss$ be a loss of the form
%
%
\begin{equation}
\label{eqnmargin-surrloss} \surrloss(\alpha, \preflabel) = \sum
_{i,j:\preflabel_{ij} > 0} \phi \bigl(\alpha_i - \alpha_j
- h(\preflabel_{ij}) \bigr).
\end{equation}
If $\phi$ is convex, then even in the
low-noise setting of Definition~\ref{defrev-triangle} the loss
$\surrloss$
is not structure-consistent.
\end{theorem}

\subsection{Achieving Fisher consistency}
\label{secachieve-consistency}

Although Section~\ref{secconsistency-hard} suggests an inherent
difficulty in
the development of tractable losses for ranking, tractable Fisher consistency
is in fact achievable if one has access to \emph{complete} preference data.
We review a few of the known results here, showing how they follow from the
Fisher consistency guarantees in Section~\ref{secgeneral-consistency}, and
derive some new Fisher consistency guarantees for the complete data setting
(we defer all proofs to
the supplementary material \cite{DuchiMaJo13supp}).
These results may appear to be of limited practical value, since complete
preference judgments are typically unavailable or untrustworthy, but,
as we
show in Sections~\ref{secu-statistics} and \ref
{secrank-aggregation}, they
can be combined with aggregation strategies to yield procedures that
are both
practical and come with consistency guarantees.

We first define the normalized discounted cumulative gain (NDCG) family
of complete data losses. Such losses are common in applications like web
search, since they penalize ranking errors at the top of a ranked list more
heavily than errors farther down the list. Let $\structure
\in\structurespace\subseteq\R^\numresults$ be a vector of
relevance scores
and $\alpha\in\R^\numresults$ be a vector of predicted scores.
Define $\pi_\alpha$ to be the permutation associated
with $\alpha$, so that $\pi_\alpha(j)$ is the rank of item $j$ in
the ordering
induced by $\alpha$.
Following Ravikumar et~al.~\cite{RavikumarTeYa11}, a
general class of NDCG loss functions can be defined as follows:
%
%
\begin{equation}
\label{eqnndcg} L(\alpha, \structure) = 1 - \frac{1}{Z(\structure)} \sum
_{j = 1}^{\numresults} \frac{G(\structure_j)}{F(\pi_\alpha(j))}, \qquad Z(\structure) = \max
_{\alpha'} \sum_{j=1}^m
\frac{G(\structure
_j)}{F(\pi_{\alpha'}(j))},
\end{equation}
where $G$ and $F$ are functions monotonically increasing in their arguments.
By inspection, $L \in[0, 1]$, and we remark that the standard NDCG
criterion \cite{JarvelinKe04} uses $G(\structure_j) = 2^{\structure
_j} - 1$
and $F(j) = \log(1 + j)$. The ``precision at $k$''
loss \cite{ManningRaSc08} can also be written in the form \eqref{eqnndcg},
where $G(\structure_j) = \structure_j$ (assuming that $\structure_j
\ge0$) and $F(j) = 1$ for $j \le k$ and $F(j) = +\infty$ otherwise, which
measures the relevance of the top $k$ items given by the vector
$\alpha$. This form generalizes standard forms of precision, which
assume $\structure_j \in\{0, 1\}$.

To analyze the consistency of surrogate losses for the NDCG
family \eqref{eqnndcg}, we first compute the loss $\condloss(\alpha,
\limitlaw)$ and then state a corollary to
Proposition~\ref{propositionall-consistency}. Observe that
for any $\limitlaw\in\measures(\structurespace)$,
\[
\condloss(\alpha, \limitlaw) = 1 - \sum_{j=1}^\numresults
\frac{1}{F(\pi_\alpha(j))} \int\frac{G(\structure_j)}{Z(\structure)} \,d\limitlaw(\structure).
\]
Since the function $F$ is increasing in its argument, minimizing
$\condloss(\alpha, \limitlaw)$ corresponds to choosing any vector
$\alpha$
whose values $\alpha_j$ obey the same order as the $\numresults$
points $\int
G(\structure_j) / Z(\structure) \,d\limitlaw(\structure)$. In
particular, the
range of $\condloss$ is finite for any $\limitlaw$ since it depends
only on
the permutation induced by $\alpha$, so we have
Corollary~\ref{corollaryndcg-consistency}.
%
\begin{corollary}
\label{corollaryndcg-consistency}
Define the set
%
%
\begin{equation}
\label{eqnndcg-optimal-set} A(\limitlaw) = \biggl\{\alpha\in\R^\numresults
\mid \alpha_j > \alpha_l \mbox{ when } \int
\frac{G(\structure_j)}{Z(\structure)} \,d\limitlaw(\structure) > \int\frac{G(\structure_l)}{Z(\structure)} \,d\limitlaw(
\structure) \biggr\}.
\end{equation}
A surrogate loss $\surrloss$ is Fisher-consistent for the NDCG
family \eqref{eqnndcg} if and only if for all $\limitlaw\in
\measures(\structurespace)$,
\[
\inf_\alpha \Bigl\{\condsurrloss(\alpha, \limitlaw) - \inf
_{\alpha'} \condsurrloss\bigl(\alpha', \limitlaw\bigr)
\mid\alpha\notin A(\limitlaw) \Bigr\} > 0.
\]
\end{corollary}

Corollary~\ref{corollaryndcg-consistency} recovers the main flavor of the
consistency results in the papers of Ravikumar et~al. \cite{RavikumarTeYa11}
and Buffoni et~al. \cite{BuffoniCaGaUs11}. The surrogate $\surrloss$ is consistent
if and
only if it preserves the order of the integrated terms $\int
G(\structure_j) /
Z(\structure) \,d\limitlaw(\structure)$: any sequence $\alpha_n$
tending to the infimum of $\condsurrloss(\alpha, \limitlaw)$ must satisfy
$\alpha_n \in A(\limitlaw)$ for large enough $n$. Zhang \cite{Zhang04a} presents
several examples of such losses; as a corollary to his
Theorem~5 (also noted by \cite{BuffoniCaGaUs11}), the loss
\[
\surrloss(\alpha, \structure) \defeq\sum_{j=1}^\numresults
\frac{G(\structure_j)}{Z(\structure)} \sum_{l = 1}^\numresults\phi(
\alpha_l - \alpha_j)
\]
is convex and structure-consistent (in the sense of
Definition~\ref{definitionstructure-consistent}) whenever $\phi\dvtx  \R
\rightarrow\R_+$ is nonincreasing, differentiable and satisfies
$\phi'(0) <
0$. The papers \cite{RavikumarTeYa11,BuffoniCaGaUs11} contain more examples
and a deeper study of NDCG losses. To extend
Corollary~\ref{corollaryndcg-consistency} to a uniform result, we
note that
if $G(\structure_j) > 0$ for all $j$ and $\structurespace$ is
compact, then
$\surrloss$ is 0-coercive over the set $\{\alpha_1 = 0\}$,\setcounter{footnote}{2}\footnote
{The loss
is invariant to linear shifts by the ones vector $\onevec$, so we may
arbitrarily set a value for $\alpha_1$.} whence
Theorem~\ref{theoremstructure-consistency} implies that structure consistency
coincides with uniform consistency.

Another family of loss functions is based on a cascade model of user
behavior~\cite{ChapelleMeZhGr09}. These losses model dependency among items
or results by assuming that a user scans an ordered list of results
from top
to bottom and selects the first satisfactory result. Here, satisfaction is
determined independently at each position. Let $\pi_\alpha^{-1}(i)$ denote
the index of item that $\alpha\in\R^\numresults$ assigns
to rank $i$.
The form of such expected reciprocal
rank (ERR) losses is
%
%
\begin{equation}
\label{eqnerr} L(\alpha, \structure) = 1 - \sum_{i = 1}^\numresults
\frac{1}{F(i)} G(\structure_{\pi_\alpha^{-1}(i)}) \prod_{j=1}^{i - 1}
\bigl(1 - G(\structure_{\pi_\alpha^{-1}(j)})\bigr),
\end{equation}
where $G \dvtx \R\rightarrow[0, 1]$ is a nondecreasing function that indicates
the prior probability that a result with score $\structure_j$ is
selected, and
$F \dvtx \N\rightarrow [{1}, {\infty} )$ is an increasing
function that
more heavily weights the first items. The ERR family also satisfies $L
\in[0,
1]$, and empirically correlates well with user satisfaction in ranking
tasks~\cite{ChapelleMeZhGr09}.

Computing the expected conditional risk $\condloss(\alpha, \limitlaw
)$ for
general $\limitlaw\in\measures(\structurespace)$ is difficult, but
we can
compute it when $\limitlaw$ is a product measure over $\structure_1,
\ldots,
\structure_\numresults$. Indeed, in this case, we have
\begin{eqnarray*}
\condloss(\alpha, \limitlaw) & = &1 - \sum_{i = 1}^\numresults
\frac{1}{F(i)} \int G(\structure_{\pi_\alpha^{-1}(i)}) \prod
_{j = 1}^{i-1} \bigl(1 - G(\structure_{\pi_\alpha^{-1}(j)})
\bigr) \,d\limitlaw(\structure)
\\
& = &1 - \sum_{i = 1}^\numresults\frac{1}{F(i)}
\E_\limitlaw\bigl[G(\structure_{\pi_\alpha^{-1}(i)})\bigr] \prod
_{j=1}^{i-1} \bigl(1 - \E_\limitlaw\bigl[G(
\structure_{\pi_\alpha^{-1}(j)})\bigr] \bigr).
\end{eqnarray*}
When one believes that the values $G(\structure_i)$ represent the a priori
relevance of the result $i$, this independence assumption is not unreasonable,
and indeed, in Section~\ref{secrank-aggregation} we provide examples
in which
it holds. Regardless, we see that $\condloss(\alpha, \limitlaw)$
depends only
on the permutation $\pi_\alpha$, and we can compute the minimizers of
the conditional risk for the ERR family \eqref{eqnerr} using the following
lemma, with proof provided in
the supplementary material \cite{DuchiMaJo13supp}.
%
\begin{lemma}
\label{lemmaerr-minimizers}
Let $p_i = \E_\limitlaw[G(\structure_i)]$. The permutation $\pi$ minimizing
$\condloss(\alpha, \limitlaw)$ is in decreasing order of the $p_i$.
\end{lemma}

Lemma~\ref{lemmaerr-minimizers} shows that an order-preserving
property is
necessary and sufficient for the Fisher-consistency of a surrogate
$\surrloss$
for the ERR family \eqref{eqnerr}, as it was for the NDCG
family \eqref{eqnndcg}. To see this, we apply a variant of
Corollary~\ref{corollaryndcg-consistency} where $A(\limitlaw)$ as
defined in
equation \eqref{eqnndcg-optimal-set} is replaced with the set
\[
A(\limitlaw) = \biggl\{\alpha\in\R^\numresults \mid\alpha_j >
\alpha_l \mbox{ whenever } \int G(\structure_j) \,d
\limitlaw(\structure) > \int G(\structure_l) \,d\limitlaw(\structure)
\biggr\}.
\]
Theorem~5 of \cite{Zhang04a} implies that
$\surrloss(\alpha, \structure) = \sum_{j=1}^\numresults
G(\structure_j)
\sum_{l = 1}^\numresults\phi(\alpha_l - \alpha_j)$
is a consistent surrogate when $\phi$ is convex,
differentiable and nonincreasing with $\phi'(0) < 0$.
Theorem~\ref{theoremstructure-consistency} also yields an equivalence
between structure and uniform consistency under suitable conditions
on $\structurespace$.

Before concluding this section, we make a final remark, which has
bearing on
the aggregation strategies we discuss in Section~\ref{secrank-aggregation}.
We have assumed that the structure spaces $\structurespace$ for the
NDCG \eqref{eqnndcg} and ERR \eqref{eqnerr} loss families consist of
real-valued relevance scores. This is certainly not necessary. In some
situations, it may be more beneficial to think of $\structure
\in\structurespace$ as simply an ordered list of the results or as a directed
acyclic graph over $\{1, \ldots, \numresults\}$. We can then apply a
transformation $r \dvtx \structurespace\rightarrow\R^\numresults$ to get
relevance scores, using $r(\structure)$ in place of $\structure$ in the
losses \eqref{eqnndcg} and \eqref{eqnerr}.
This has the advantage of causing
$\structurespace$ to be finite, so Theorem~\ref{theoremstructure-consistency}
applies, and there exists a nondecreasing function $\zeta$ with
$\zeta(0) =
0$ such that for any distribution and any measurable~$f$,
\[
\risk(f) - \risk^* \le\zeta\bigl(\risk_\surrloss(f) -
\risk_\surrloss^*\bigr).
\]

\section{Uniform laws and asymptotic consistency}
\label{secu-statistics}

In Section~\ref{secconsistency}, we gave examples of losses based on readily
available pairwise data but for which Fisher-consistent tractable surrogates
do not exist. The existence of Fisher-consistent tractable surrogates for
other forms of data, as in Section~\ref{secachieve-consistency}, suggests
that aggregation of pairwise and partial data into more complete data
structures, such as lists or scores, makes the problem easier. However,
it is
not obvious how to design statistical procedures based on
aggregation. In this section, we formally define a class of suitable
estimators that permit us to take advantage of the weak convergence of
Assumption~\ref{assumptionlimit-structures} and show that uniform
laws of
large numbers hold for our surrogate losses. This means that we can indeed
asymptotically minimize the risk~\eqref{eqnrisk-def} as desired.

Our aim is to develop an empirical analogue of the population surrogate
risk \eqref{eqnsurrogate-risk} that converges uniformly to the population
risk under minimal assumptions on the loss $\surrloss$ and structure function
$\structure$. Given a dataset $\{(Q_i, \preflabel_i)\}_{i=1}^n$ with $(Q_i,
\preflabel_i) \in\queryspace\times\preflabelspace$, we begin by defining,
for each query $q$, the batch of data belonging to the query, $\batch
(q) = \{i
\in\{1,\ldots,n\} \mid Q_i = q\}$, and the empirical count of the
number of
items in the batch, $\what{n}_q = |\batch(q)|$. As a first attempt at
developing
an empirical objective, we might consider an empirical surrogate risk
based on
complete aggregation over the batch of data belonging to each query:
%
%
\begin{equation}
\label{eqnfull-aggregation-risk} \frac{1}{n} \sum
_q \what{n}_q \surrloss \bigl(f(q), \structure
\bigl(\bigl\{ \preflabel_{i_j} 
\mid i_j \in
\batch(q)\bigr\}\bigr) \bigr).
\end{equation}
While we would expect this risk to converge uniformly when $\surrloss$
is a
sufficiently smooth function of its structure argument, the analysis of the
complete aggregation risk \eqref{eqnfull-aggregation-risk} requires overly
detailed knowledge of the surrogate $\surrloss$ and the structure function
$\structure$.

To develop a more broadly applicable statistical procedure, we instead
consider an empirical surrogate based on $U$-statistics. By trading off
the nearness of an order-$k$ $U$-statistic to an i.i.d. sample and the
nearness of the limiting structure distribution $\limitlaw_q$ to a structure
$s(\preflabel_1, \ldots, \preflabel_k)$ aggregated over $k$ draws,
we can
obtain consistency under mild assumptions on $\surrloss$ and
$\structure$. More specifically, for each query~$q$, we consider the
surrogate loss
%
%
\begin{equation}
\label{eqnu-statistic-loss} \pmatrix{\what{n}_q
\cr
k}^{-1} \mathop{\sum_{i_1 < \cdots< i_k,}}_{i_j \in\batch(q)} \surrloss \bigl(f(q),
\structure(\preflabel_{i_1}, \ldots, \preflabel_{i_k}) \bigr).
\end{equation}
When $\what{n}_q < k$, we adopt the convention ${\what{n}_q
\choose k} = 1$, and
the above sum becomes the single term $\surrloss(f(q),
\structure(\{\preflabel_{i_j} \mid i_j \in\batch(q)\}))$ as in the
expression \eqref{eqnfull-aggregation-risk}. Hence, our $U$-statistic loss
recovers the complete aggregation loss \eqref
{eqnfull-aggregation-risk} when
$k = \infty$.

An alternative formulation to loss \eqref{eqnu-statistic-loss} might consist
of $ \lceil{|\batch(q)| / k} \rceil$ aggregation terms
per query, with each
query-preference pair appearing in a single term. However, the
instability of
such a strategy is high: a change in the ordering of the data or a
substitution of queries could have a large effect on the final
estimator. The
$U$-statistic \eqref{eqnu-statistic-loss} grants robustness to such
perturbations in the data. Moreover, by choosing the right rate of increase
of the aggregation order $k$ as a function of $n$, we obtain consistent
procedures for a broad class of surrogates $\surrloss$ and structures
$\structure$.

We associate with the surrogate loss \eqref{eqnu-statistic-loss} a surrogate
empirical risk that weights each query by its empirical probability of
appearance:
%
%
\begin{equation}
\label{eqnu-statistic-emprisk}\qquad \emprisk_{\surrloss,n}(f) \defeq
\frac{1}{n} \sum_q \what{n}_q
\pmatrix{\what{n}_q
\cr
k}^{-1}\mathop{\sum_{i_1 < \cdots< i_k,}}_{i_j \in\batch(q)} \surrloss \bigl(f(q), \structure(\preflabel_{i_1},
\ldots, \preflabel_{i_k}) \bigr).
\end{equation}
Let $\P_n$ denote the probability distribution of the queries given
that the
dataset is of size $n$. Then by iteration of expectation and Fubini's theorem,
the surrogate risk~\eqref{eqnu-statistic-emprisk} is an unbiased
estimate of
the population quantity
%
%
\begin{equation}
\label{eqnu-statistic-risk}\quad  \risk_{\surrloss,n}(f) \defeq \sum
_q \Biggl[ \sum_{l = 1}^n
l \P_n(\what{n}_q = l) \E \bigl[\surrloss \bigl(f(Q),
\structure(\preflabel_1, \ldots, \preflabel_{l \wedge k}) \bigr)
\Bigm| Q = q \bigr] \Biggr].
\end{equation}

It remains to establish a uniform law of large numbers guaranteeing the
convergence of the empirical risk \eqref{eqnu-statistic-emprisk} to
the target
population risk \eqref{eqnsurrogate-risk}. Under suitable
conditions such as those of Section~\ref{secconsistency}, this
ensures the
asymptotic consistency of computationally tractable statistical
procedures. Hereafter, we assume that we have a nondecreasing sequence of
function classes $\funclass_n$, where any $f \in\funclass_n$ is a scoring
function for queries, mapping $f \dvtx \queryspace\rightarrow\R
^\numresults$ and
giving scores to the (at most $\numresults$) results for each\vadjust{\goodbreak} query $q
\in\queryspace$. Our goal is to give sufficient conditions for
the convergence in probability
%
%
\begin{equation}
\label{eqnuniform-convergence} \sup_{f \in\funclass_n} \bigl\llvert
\emprisk_{\surrloss, n}(f) - \risk_\surrloss(f) \bigr\rrvert \cp0\qquad \mbox{as }
n \rightarrow\infty.
\end{equation}

While we do not provide fully general conditions under which the
convergence~\eqref{eqnuniform-convergence} occurs, we provide representative,
checkable conditions sufficient for convergence. At a high
level, to establish \eqref{eqnuniform-convergence}, we control
the uniform difference between the expectations $\risk_{\surrloss, n}(f)$
and $\risk_\surrloss(f)$ and bound the distance between the empirical risk
$\emprisk_{\surrloss,n}$ and its expectation $\risk_{\surrloss, n}$
via covering number arguments.
We now specify assumptions under
which our results hold, deferring all proofs to
the supplementary material~\cite{DuchiMaJo13supp}.

Without loss of generality, we assume that $\queryprob_q$, the true
probability of seeing the query $q$, is nonincreasing in the query
index $q$.
First, we describe the tails of the query distribution:
%
\renewcommand{\theassumption}{\Alph{assumption}}
\setcounter{assumption}{3}
\begin{assumption}
\label{assumptionpower-law}
There exist constants $\queryexp> 0$ and $K_1 > 0$ such that
$\queryprob_q \le K_1 q^{-\queryexp- 1}$ for all $q$. That is,
$\queryprob_q = \order(q^{-\queryexp- 1})$.
\end{assumption}

Infinite sets of queries $\queryspace$ are reasonable, since search
engines, for example, receive a large volume of entirely new queries
each day.
Our arguments also apply when $\queryspace$ is finite, in which case
we can take $\queryexp\uparrow\infty$.

Our second main assumption concerns the behavior of
the surrogate loss $\surrloss$ over the function class $\funclass_n$,
which we
assume is contained in a normed space with norm~$\llVert {\cdot}\rrVert $.
%
\begin{assumption}[(Bounded Lipschitz losses)]
\label{assumptionbounded-lipschitz}
The surrogate loss function $\surrloss$ is bounded and Lipschitz continuous
over $\funclass_n$: for any $\structure\in\structurespace$,
any $f, f_1, f_2 \in\funclass_n$, and any $q \in\queryspace$,
there exist constants $\surrbound_n$ and $\lipbound_n < \infty$ such that
\[
0 \le\surrloss\bigl(f(q), \structure\bigr) \le\surrbound_n
\]
and
\[
\bigl\llvert \surrloss\bigl(f_1(q), \structure\bigr) - \surrloss
\bigl(f_2(q), \structure\bigr)\bigr\rrvert \le\lipbound_n
\llVert {f_1 - f_2}\rrVert .
\]
\end{assumption}

This assumption is satisfied whenever $\surrloss(\cdot, \structure)$
is convex
and $\funclass_n$ is compact [and contained in the interior of the
domain of
$\surrloss(\cdot, \structure)$] \cite{HiriartUrrutyLe96ab}. Our final
assumption gives control over the sizes of the function classes
$\funclass_n$
as measured by their covering numbers. (The $\varepsilon$-covering number of
$\funclass$ is the smallest $N$ for which there are $f^i$, $i \le N$, such
that $\min_i \|{f^i - f}\| \le\varepsilon$ for any $f \in\funclass$.)
%
\begin{assumption}
\label{assumptioncovering-number}
For all $\varepsilon> 0$, $\funclass_n$ has $\varepsilon$-covering
number $N(\varepsilon, n) < \infty$.
\end{assumption}


With these assumptions in place, we give a few representative
conditions that
enable us to guarantee uniform convergence \eqref{eqnuniform-convergence}.\vadjust{\goodbreak}
Roughly, these conditions control the interaction between the size
of the function classes $\funclass_n$ and the order $k$ of aggregation used
with $n$ data points. To that end, we let the aggregation order $k_n$ grow
with $n$. In stating the conditions, we make use of the shorthand
$\E_q[\surrloss(f(q), \structure(\preflabel_{1:k}))]$ for
$\E[\surrloss(f(Q), \structure(\preflabel_1, \ldots, \preflabel_k))
\mid Q = q]$.
%
\renewcommand{\thecondition}{\Roman{condition}}
\setcounter{condition}{0}
\begin{condition}
\label{conditionpolynomial-in-k}
There exist a $\rho> 0$ and constant $C$ such that for all
$q \in\queryspace$, $n \in\N$, $k \in\N$, and $f \in\funclass_n$,
\[
\Bigl\llvert \E_q\bigl[\surrloss\bigl(f(q), \structure(
\preflabel_1, \ldots, \preflabel_k)\bigr)\bigr] - \lim
_{k'} \E_q\bigl[\surrloss\bigl(f(q), \structure(
\preflabel_1, \ldots, \preflabel_{k'})\bigr)\bigr] \Bigr
\rrvert \le C \surrbound_n k^{-\rho}.
\]
Additionally, the sequences $\surrbound_n$ and $k_n$
satisfy $\surrbound_n = o(k_n^\rho)$.
\end{condition}

This condition is not unreasonable; when $\surrloss$ and $\structure$ are
suitably continuous, we expect $\rho\ge\frac{1}{2}$. We also
consider an
alternative covering number condition.

\renewcommand{\thecondition}{\Roman{condition}$^{\prime}$}
\setcounter{condition}{0}
\begin{condition}
\label{conditioncovering-expectations}
Sequences $\{\varepsilon_n\} \subset\R_+$ and $\{k_n\} \subset\N$ and an
$\varepsilon_n$-cover $\funclass_n^1, \ldots,  \funclass_n^{N(\varepsilon
_n, n)}$
of $\funclass_n$ can be chosen such that
\[
\max_{i \in\{1,\dots,N(\varepsilon_n, n)\}} \inf_{f \in\funclass_n^i} \biggl|
\risk_\surrloss(f) - \sum_q
\queryprob_q \E_q\bigl[\surrloss\bigl(f(q), \structure(
\preflabel_1, \ldots, \preflabel_{k_n})\bigr)\bigr]\biggr | + 2
\lipbound_n \varepsilon_n \rightarrow0.
\]
\end{condition}
%

Condition \ref{conditioncovering-expectations} is weaker than
Condition \ref{conditionpolynomial-in-k}, since it does not require uniform
convergence over $q \in\queryspace$. If the function class $\funclass
$ is
fixed for all $n$, then the weak convergence of $\structure(\preflabel_1,
\ldots, \preflabel_k)$ as in Assumption~\ref{assumptionlimit-structures}
guarantees Condition~\ref{conditioncovering-expectations}, since
$N(\varepsilon,
n) = N(\varepsilon, n') < \infty$, and we may take $\varepsilon$
arbitrarily small.
We require one additional condition, which relates the growth of $k_n$,
$\surrbound_n$, and the function classes $\funclass_n$ more directly.
%
\renewcommand{\thecondition}{\Roman{condition}}
\setcounter{condition}{1}
\begin{condition}
\label{conditionall-sequence-growths}
The sequences $k_n$ and $\surrbound_n$ satisfy $k_n \surrbound
_n^{{(1 + \queryexp)}/{\queryexp}} = o(n)$. Additionally, for any fixed
$\varepsilon
> 0$, the sequences satisfy
\[
k_n \surrbound_n \biggl[\log N \biggl(
\frac{\varepsilon}{4 \lipbound_n}, n \biggr) \biggr]^{{1}/ {2}}= o(\sqrt{n}).
\]
\end{condition}

By inspection, Condition \ref{conditionall-sequence-growths} is
satisfied for
any $k_n = o(\sqrt{n})$ if the function classes $\funclass_n$ are
fixed for
all $n$. Similarly, if for all $k \ge k_0$, $\structure(\preflabel_1,
\ldots,
\preflabel_k) = \structure(\preflabel_1, \ldots, \preflabel
_{k_0})$, so
$\structure$ depends only on its first $k_0$ arguments,
Condition \ref{conditionall-sequence-growths} holds whenever
$\max\{\surrbound_n^{(1 + \queryexp)/\queryexp}, \surrbound_n^2
\log
N(\varepsilon/4\lipbound_n, n)\} = o(n)$. If the function classes
$\funclass_n$
consist of linear functionals represented by vectors $\theta\in\R
^{d_n}$ in
a ball of some finite radius, then $\log N(\varepsilon, n) \approx d_n
\log
\varepsilon^{-1}$, which means that
Condition~\ref{conditionall-sequence-growths} roughly requires $k_n
\sqrt{d_n
/ n} \rightarrow0$ as $n \rightarrow\infty$. Modulo the factor
$k_n$, this
condition is familiar from its necessity in the convergence of parametric
statistical problems.

The conditions in place, we come to our main result on the convergence of
our $U$-statistic-based empirical loss minimization procedures.
%
\begin{theorem}
\label{theoremconvergence}
Assume Condition \ref{conditionpolynomial-in-k}
or \ref{conditioncovering-expectations} and additionally assume
the growth Condition \ref{conditionall-sequence-growths}. Under
Assumptions \ref{assumptionpower-law}, \ref{assumptionbounded-lipschitz}
and \ref{assumptioncovering-number},
\[
\sup_{f \in\funclass_n} \bigl\llvert \emprisk_{\surrloss, n}(f) -
\risk_\surrloss(f)\bigr\rrvert \cp0 \qquad\mbox{as } n \rightarrow\infty.
\]
\end{theorem}

We remark in passing that if Condition \ref{conditionall-sequence-growths}
holds, with the change that the $o(\sqrt{n})$ bound is replaced by
$O(n^{\rho})$
for some $\rho< \frac{1}{2}$, the conclusion of
Theorem~\ref{theoremconvergence} can be strengthened to both
convergence almost surely and in expectation.

By inspection, Theorem~\ref{theoremconvergence} provides our desired
convergence guarantee \eqref{eqnuniform-convergence}. By combining the
Fisher-consistent loss families outlined in
Section~\ref{secachieve-consistency} with the consistency guarantees provided
by Theorem~\ref{theoremconvergence}, it is thus possible to design
statistical procedures that are both computationally tractable---minimizing
only convex risks---and asymptotically consistent.

\section{Rank aggregation strategies}
\label{secrank-aggregation}

In this section, we give several examples of practical strategies for
aggregating disparate user preferences under our framework. Motivated
by the
statistical advantages of complete preference data highlighted in
Section~\ref{secachieve-consistency}, we first present strategies
for constructing complete vectors of relevance scores from pairwise preference
data. We then discuss a model for the selection or ``click'' data that arises
in web search and information retrieval and show that maximum likelihood
estimation under this model allows for consistent ranking. We conclude this
section with a brief overview of structured aggregation strategies.

\subsection{Recovering scores from pairwise preferences}

Here we treat partial preference observations as noisy evidence of an
underlying complete ranking and attempt to achieve consistency with
respect to
a complete preference data loss. We consider three methods that take as input
pairwise preferences and output a relevance score vector $\structure
\in
\R^\numresults$. Such procedures fit naturally into our
ranking-with-aggregation
framework: the results in Section~\ref{secachieve-consistency}
and Section~\ref{secu-statistics} show that a Fisher-consistent loss is
consistent for the limiting distribution of the scores $\structure$ produced
by the aggregation procedure. Thus, it is the responsibility of the
statistician---the designer of an aggregation procedure---to determine whether
the scores accurately reflect the judgments of the population. We
present our
first example in some detail to show how aggregation of pairwise
judgments can
lead to consistency in our framework and follow with brief descriptions of
alternate aggregation strategies. For an introduction to the design of
aggregation strategies for pairwise data, see Tsukida and Gupta~\cite{TsukidaGu11} as
well as
the book by David \cite{David69}.

\textit{Thurstone--Mosteller least squares and skew-symmetric scoring}.
The first aggregation strategy constructs a relevance score vector
$\structure$ in two phases. First, it aggregates a sequence of observed
preference judgments $\preflabel_i \in\preflabelspace$, provided in any
form, into a skew-symmetric matrix $A \in\R^{\numresults\times
\numresults}$ satisfying $A = -A^\top$. Each entry $A_{ij}$ encodes the
extent to which item $i$ is preferred to item $j$. Given such a
skew-symmetric matrix, Thurstone and Mosteller \cite{Mosteller51} recommend
deriving a score vector $\structure$ such that $\structure_i -
\structure_j
\approx A_{ij}$. In practice, one may not observe preference
information for
every pair of results, so we define a masking matrix $\Omega\in
\{0,1\}^{\numresults\times\numresults}$ with $\Omega= \Omega^\top$,
$\Omega_{ii}=1$, and $\Omega_{ij} = 1$ if and only if preference information
has been observed for the pair $i\neq j$. Letting $\circ$ denote the Hadamard
product, a natural objective for selecting
scores (e.g., \cite{Gulliksen56})
is the least squares objective
%
%
\begin{equation}
\label{eqnskew-symmetric-score}
\qquad\mathop{\minimize}_{x \dvtx x^\top {\mathbh{1}} = 0} \frac{1}{4} \sum
_{i,j} \Omega_{ij} \bigl(A_{ij} -
(x_i - x_j)\bigr)^2 = \frac{1}{4}
\bigl\llVert {\Omega\circ\bigl(A - \bigl(\onevec x^\top- x \onevec
^\top\bigr)\bigr)}\bigr\rrVert _\mathrm{Fr}^2.
\end{equation}
The gradient of the objective \eqref{eqnskew-symmetric-score} is
\[
D_\Omega x - (\Omega\circ A) \onevec- \Omega x \qquad\mbox{where }
D_\Omega\defeq\diag(\Omega\onevec).
\]
Setting $\structure= (D_\Omega- \Omega)^\pinv(\Omega
\circ A)\onevec$ yields the solution to the minimization
problem~\eqref{eqnskew-symmetric-score}, since
$D_\Omega- \Omega$ is an unnormalized graph
Laplacian matrix \cite{Chung98}, and therefore
$\onevec^\top s = \onevec^\top(D_\Omega- \Omega)^\pinv(\Omega
\circ A)\onevec= 0$.

If $\Omega= \onevec\onevec^\top$, so that all pairwise preferences are
observed, then the eigenvalue decomposition of $D_\Omega- \Omega=
\numresults I - \onevec\onevec^\top$ can be computed explicitly as
$V \Sigma
V^\top$, where $V$ is any orthonormal matrix whose first column is $1 /
\sqrt{\numresults}$, and $\Sigma$ is a diagonal matrix with entries
$0$ (once)
and $\numresults$ repeated $\numresults- 1$ times. Thus, letting
$\structure_A$ and $\structure_B$ denote solutions to the minimization
problem \eqref{eqnskew-symmetric-score} with different skew-symmetric
matrices $A$ and $B$ and noting that $A \onevec\perp\onevec$ since
$\onevec^\top A \onevec= 0$, we have the Lipschitz continuity of the
solutions $\structure$ in $A$:
\[
\llVert {\structure_A - \structure_B}\rrVert
_2^2 = \bigl\llVert {\bigl(\numresults I - \onevec
\onevec^\top\bigr)^\pinv(A - B)\onevec}\bigr\rrVert
_2^2 = \frac{1}{\numresults^2} \bigl\llVert {(A - B)\onevec}
\bigr\rrVert _2^2 \le\frac{1}{\numresults} \llvert\!
\llvert \!\llvert {A - B}\rrvert \!\rrvert\! \rrvert _2^2.
\]
Similarly,
when $\Omega$ is fixed, the score structure $\structure$ is likewise
Lipschitz in $A$ for any norm $\llvert\!  \llvert  \!\llvert {\cdot}\rrvert\!  \rrvert\!  \rrvert $ on skew-symmetric matrices.

A variety of procedures are available for aggregating pairwise
comparison data
$\preflabel_i \in\preflabelspace$ into a skew-symmetric matrix $A$.
One example,
the Bradley--Terry--Luce (BTL) model \cite{BradleyTe52}, is based
upon empirical log-odds ratios. Specifically, assume that $\preflabel_i
\in\preflabelspace$ are pairwise comparisons of the form $j \succ l$, meaning
item $j$ is preferred to item $l$. Then we can set
\[
A_{jl} = \log\frac{\what{\P}(j \succ l) + c}{\what{\P}(j \prec l)
+ c} \qquad\mbox{for observed pairs } j, l,
\]
where $\what{\P}$ denotes the empirical distribution over $\{
\preflabel_1,
\ldots, \preflabel_k\}$ and $c > 0$ is a smoothing parameter.

Since the proposed structure $\structure$ is a continuous function of the
skew-symmetric matrix $A$, the limiting distribution $\limitlaw$ is a point
mass whenever $A$ converges almost surely, as it does in the BTL model.
If aggregation is carried out using only a finite number of
preferences rather than letting $k$ approach $\infty$ with $n$, then
$\limitlaw$ converges to a nondegenerate distribution.
Theorem~\ref{theoremstructure-consistency} grants uniform consistency
since the score space $\structurespace$ is finite.

\textit{Borda count and budgeted aggregation.}
The Borda count \cite{Borda1781} provides a computationally efficient method
for computing scores from election results. In a general election setting,
the procedure counts the number of times that a particular item was
rated as
the best, second best, and so on. Given a skew-symmetric matrix $A$
representing the outcomes of elections, the Borda count assigns the scores
$\structure= A\onevec$. As above, a skew-symmetric matrix $A$ can be
constructed from input preferences $\{\preflabel_1, \ldots,
\preflabel_k\}$,
and the choice of this first-level aggregation can greatly affect the
resulting rankings. Ammar and Shah \cite{AmmarSh11} suggest that if one has limited
computational budget and only pairwise preference information then one should
assign to item $j$ the score
\[
\structure_j = \frac{1}{\numresults- 1} \sum
_{l \neq j} \what{\P}(j \succ l), \label{eqnammar-score}
\]
which estimates the probability of winning
an election against an opponent chosen uniformly.
This is equivalent to the Borda count
when we choose $A_{jl} = \what{\P}(j \succ l) - \what{\P}(j \prec l)$
as the entries in the skew-symmetric aggregate $A$.

\textit{Principal eigenvector method.}
Saaty \cite{Saaty03} describes the principal eigenvector method, which
begins by
forming a reciprocal matrix $A \in\R^{\numresults\times\numresults
}$, with
positive entries $A_{ij} = (A_{ji})^{-1}$, from pairwise comparison judgments.
Here $A_{ij}$ encodes a multiplicative preference for item $i$ over
item $j$;
the idea is that ratios preserve preference strength \cite{Saaty03}. To
generate $A$, one may use, for example, smoothed empirical ratios
$A_{jl} =
\frac{\what{\P}(j \succ l) + c}{\what{\P}(j \prec l) + c}$.
Saaty recommends finding a vector $\structure$ so that
$\structure_i / \structure_j \approx A_{ij}$, suggesting using the Perron
vector of the matrix, that is, the first eigenvector of $A$.

\subsection{Cascade models for selection data}
\label{seccascade-models}
\newcommand{\result}{l}

Cascade models \cite{CraswellZoTaRa08,ChapelleMeZhGr09} explain
the behavior of a user presented with an ordered list of items, for example
from a web search. In a cascade model, a user considers results in the
presented order and selects the first to satisfy him or her.
The model assumes the result $\result$ satisfies a user
with probability $p_\result$, independently of previous items in the list.
It is natural to express a variety of ranking losses, including the expected
reciprocal rank (ERR) family \eqref{eqnerr}, as expected disutility
under a
cascade model, but computation and optimization of these losses require
knowledge of the satisfaction probabilities $p_\result$. When the
satisfaction probabilities are unknown, Chapelle et~al. \cite{ChapelleMeZhGr09} recommend
plugging in those values $p_\result$ that maximize the likelihood of observed
click data. Here we show that risk consistency for the ERR family is
straightforward
to characterize when scores are estimated via maximum likelihood.

To this end, fix a query $q$, and let each affiliated preference judgment
$\preflabel_i$ consist of a triple $(\numresults_i, \pi_i, c_i)$, where
$\numresults_i$ is the number of results presented to the user, $\pi
_i$ is the
order of the presented results, which maps positions $\{1,\ldots,
\numresults_i\}$ to the full result set $\{1,\ldots,\numresults\}$,
and $c_i
\in\{1,\ldots,\numresults_i + 1\}$ is the position clicked on by the
user ($\numresults_i + 1$ if the user chooses nothing). The likelihood
$g$ of
an i.i.d. sequence $\{\preflabel_1,\ldots, \preflabel_k\}$ under a
cascade model $p$ is
\[
g\bigl(p, \{\preflabel_1, \ldots, \preflabel_k\}\bigr) =
\prod_{i=1}^k p_{\pi_i(c_i)}^{1 ({c_i\leq m_i} )}
\prod_{j=1}^{c_i - 1} (1 - p_{\pi_i(j)}),
\]
and the maximum likelihood estimator of the satisfaction probabilities
has the
closed form
\[
\what{p}_\result(\preflabel_1,\ldots, \preflabel_k)
= \frac{\sum_{i=1}^k 1 ({\pi_i(c_i) = \result} )}{\sum_{i=1}^k \sum_{j=1}^{c_i}1 ({\pi_i(j)=\result} )},\qquad l = 1, \ldots, \numresults.
\]

To incorporate this maximum likelihood aggregation procedure into our
framework, we define the structure function $\structure$ to be the vector
\[
\structure(\preflabel_1, \ldots, \preflabel_k) \defeq
\what{p}(\preflabel_1,\ldots,\preflabel_k) \in\R
^\numresults
\]
of maximum likelihood probabilities, and we take as our loss $L$ any
member of
the ERR family \eqref{eqnerr}. The strong law of large numbers
implies the
a.s. convergence of $\what{p}$ to a vector $p \in[0,1]^\numresults
$, so that
the limiting law $\limitlaw_q(\{p\}) = 1$. Since $\limitlaw_q$ is a product
measure over $[0, 1]^\numresults$, Lemma~\ref{lemmaerr-minimizers} implies
that any $\alpha$ inducing the same ordering over results as $p$
minimizes the
conditional ERR risk $\condloss(\alpha, \limitlaw)$.
By application of Theorems \ref{theoremstructure-consistency} (or
Proposition~\ref{propositionall-consistency}) and \ref{theoremconvergence},
it is possible to asymptotically minimize the expected reciprocal rank
by aggregation.

\subsection{Structured aggregation}
\label{secstructured-aggregation}
Our framework can leverage aggregation
procedures (see, e.g., \cite{DworkKuNaSi01}) that
map input preferences into representations of combinatorial objects.
Consider the setting of Section~\ref{secconsistency-hard}, in
which each observed preference judgment $\preflabel$ is the weighted adjacency
matrix of a directed acyclic graph, our loss of interest $L$ is the edgewise
indicator loss \eqref{eqnedge-loss}, and our candidate surrogate
losses have
the form \eqref{eqnu-statistic-loss}.
Theorems \ref{theoreminconsistency-convex}
and \ref{theoreminconsistency-margin} establish that risk consistency
is not
generally attainable when $\structure(\preflabel_1, \ldots
,\preflabel_k) =
\preflabel_1$.
In certain cases, aggregation can recover consistency. Indeed,
define
\[
\structure(\preflabel_1, \ldots, \preflabel_k) \defeq
\frac{1}{k}\sum_{i=1}^k
\preflabel_i,
\]
the average of the input adjacency matrices. For an i.i.d. sequence
$\preflabel_1,\preflabel_2,\ldots$ associated with a given query
$q$, we have\vadjust{\goodbreak}
$\structure(\preflabel_1, \ldots, \preflabel_n) \cas\E(\preflabel
\mid Q=q)$
by the strong law of large numbers, and hence the asymptotic surrogate risk
\[
\risk_\surrloss(f) = \sum_q
\queryprob_q \int\surrloss\bigl(f(q),\structure\bigr) \,d\limitlaw
_q(\structure) = \sum_q
\queryprob_q \surrloss\bigl(f(q),\E(\preflabel\mid Q=q)\bigr).
\]
Recalling the conditional pairwise risk \eqref{eqnedge-condloss},
we can rewrite the risk as
\begin{eqnarray*}
\risk(f)& =& \sum_q
\queryprob_q \biggl[\sum_{i < j}
\preflabel_{ij}^{\limitlaw_q} 1 \bigl({f_i(q) \le
f_j(q)} \bigr) + \sum_{i > j}
\preflabel_{ij}^{\limitlaw_q} 1 \bigl({f_i(q) <
f_j(q)} \bigr) \biggr]
\\
& =& \sum_q
\queryprob_q \sum_{i > j} \E[
\preflabel_{ij} \mid Q = q] \\
&&{} + \sum_q
\queryprob_q \sum_{i < j} \E[
\preflabel_{ij} - \preflabel_{ji} \mid Q = q] 1
\bigl({f_i(q) \le f_j(q)} \bigr). 
\end{eqnarray*}

The discussion immediately following
Proposition~\ref{propositionall-consistency} shows that any consistent
surrogate $\surrloss$ must be bounded away from its minimum for
$\alpha\notin\break\argmin_{\alpha'}\condsurrloss(\alpha', \limitlaw)$. Since the limiting
distribution $\limitlaw$ is a point mass at some adjacency matrix
$\structure$
for each $q$, a surrogate loss $\surrloss$ is consistent if and only if
\[
\inf_{\alpha} \Bigl\{ \surrloss(\alpha, \structure) - \inf
_{\alpha'} \surrloss\bigl(\alpha', \structure\bigr)
\mid\alpha\notin\mathop{\argmin}_{\alpha'} L\bigl(\alpha', \structure
\bigr) \Bigr\} > 0.
\]
In the important special case when the difference graph $G_\limitlaw$
associated with $\E[\preflabel\mid Q=q]$
is a DAG for each query $q$ (recall
Section~\ref{seclow-noise-inconsistency}), structure consistency is
obtained if for
each $\alpha^* \in\argmin_\alpha\surrloss(\alpha, \structure)$,
$\sign(\alpha_i^* -
\alpha_j^*) = \sign(s_{ij}-s_{ji})$ for each pair of results $i, j$.
As an example, in this setting
%
%
\begin{equation}
\label{eqnpairwise-consistent} \surrloss(\alpha, \structure) \defeq \sum
_{i, j} [{\structure_{ij} -
\structure_{ji}} ]_+ \phi(\alpha_i - \alpha_j)
\end{equation}
is consistent when $\phi$ is nonincreasing,
convex, and has derivative $\phi'(0) < 0$.


The Fisher-consistent loss \eqref{eqnpairwise-consistent} is similar
to the
inconsistent losses \eqref{eqnsurr-edge-loss} considered in
Section~\ref{secconsistency-hard}, but the coefficients adjoining each
$\phi(\alpha_i - \alpha_j)$ summand exhibit a key difference. While the
inconsistent losses employ coefficients based solely on the average $i
\to j$
weight $s_{ij}$, the consistent loss coefficients are nonlinear
functions of
the edge weight differences $s_{ij}-s_{ji}$: they
are precisely the edge weights of the
difference graph $G_\limitlaw$ introduced
Section~\ref{seclow-noise-inconsistency}. Since at least one of the two
coefficients $[s_{ij}-s_{ji}]_+$ and $[s_{ji}-s_{ij}]_+$ is always
zero, the
loss \eqref{eqnpairwise-consistent} penalizes misordering either edge
$i\to
j$ or $j\to i$. This contrasts with the inconsistent surrogates of
Section~\ref{secconsistency-hard}, which
simultaneously associate nonzero convex losses with opposing edges
$i\to j$
and $j\to i$.
Note also that our argument for the consistency of the loss
\eqref{eqnpairwise-consistent} does not require
Definition~\ref{defrev-triangle}'s low-noise assumption: consistency holds
under the weaker condition that, on average, a population's preferences are
acyclic.


\section{Experimental study and implementation}
\label{secexperiments}
\newcommand{\regularizer}{\Phi}

In this section, we describe strategies for solving the convex programs that
emerge from our aggregation approach to ranking and demonstrate
the empirical utility of our proposed procedures.
We begin with a broad
description of implementation strategies and end with a presentation of
specific experiments.

\subsection{Minimizing the empirical risk}

At first glance, the empirical risk \eqref{eqnu-statistic-emprisk} appears
difficult to minimize, since the number of terms grows exponentially in
the level of aggregation $k$. Fortunately, we may leverage techniques
from the
stochastic optimization literature \cite
{NemirovskiJuLaSh09,DuchiSi09c} to
minimize the risk \eqref{eqnu-statistic-emprisk} in time linear in $k$
and independent of $n$.
Let us consider minimizing a function of the form
%
%
\begin{equation}
\label{eqnstochastic-objective} \emprisk_N(f) \defeq
\frac{1}{N} \sum_{i=1}^N \surrloss
\bigl(f, \structure^i\bigr) + \regularizer(f),
\end{equation}
where $\{\structure^i\}_{i=1}^N$ is some collection of data,
$\surrloss(\cdot,
\structure)$ is convex in its first argument, and $\regularizer$ is a convex
regularizing function (possibly zero).

Duchi and Singer \cite{DuchiSi09c}, using ideas similar to those
of Nemirovski et~al. \cite{NemirovskiJuLaSh09}, develop a specialized stochastic gradient
descent method for minimizing composite objectives of the
form \eqref{eqnstochastic-objective}. Such methods maintain
a parameter $f^t$, which is assumed to live in convex subset $\funclass
$ of a
Hilbert space with inner product $ \langle \cdot, \cdot
\rangle $, and iteratively update
$f^t$ as follows. At iteration $t$, an index $i_t \in[N]$ is chosen
uniformly at random and the gradient $\nabla_f \surrloss(f^t,
\structure^{i_t})$ is computed at $f^t$. The parameter $f$ is then
updated via
%
%
\begin{equation}
\label{eqnfobos-update} f^{t + 1} = \mathop{\argmin}_{f \in\funclass} \biggl\{
\bigl\langle f, \nabla\surrloss\bigl(f^t, \structure^{i_t}
\bigr) \bigr\rangle + \regularizer(f) + \frac{1}{2 \eta_t} \bigl\llVert {f -
f^t}\bigr\rrVert ^2 \biggr\},
\end{equation}
where $\eta_t > 0$ is an iteration-dependent stepsize and $\llVert
{\cdot}\rrVert $
denotes the Hilbert norm. The convergence guarantees of the
update \eqref{eqnfobos-update} are
well understood \cite{NemirovskiJuLaSh09,DuchiSi09c,DuchiShSiTe10}. Define
$\overline{f}^T = (1/T) \sum_{t=1}^T f^t$ to be the average parameter after
$T$ iterations. If the function $\emprisk_N$ is strongly
convex---meaning it
has at least quadratic curvature---the step-size choice $\eta_t
\propto1 / t$
gives
\[
\E \bigl[\emprisk_N\bigl(\overline{f}^T\bigr) \bigr] -
\inf_{f \in\funclass} \emprisk_N(f) = \order \biggl(
\frac{1}{T} \biggr),
\]
where the expectation is taken with respect to the indices $i_t$
chosen during each iteration of the algorithm.
In the convex case (without assuming any stronger properties than convexity),
the step-size choice $\eta_t \propto1 / \sqrt{t}$ yields
\[
\E \bigl[\emprisk_N\bigl(\overline{f}^T\bigr) \bigr] -
\inf_{f \in\funclass} \emprisk_N(f) = \order \biggl(
\frac{1}{\sqrt{T}} \biggr).
\]
These guarantees also hold with high probability \cite{NemirovskiJuLaSh09,DuchiShSiTe10}.

Neither of the convergence rates $1 / T$ or $1 / \sqrt{T}$ depends on the
number of terms $N$ in the stochastic
objective \eqref{eqnstochastic-objective}. As a consequence, we can
apply the
composite stochastic gradient method \eqref{eqnfobos-update} directly
to the
empirical risk \eqref{eqnu-statistic-emprisk}: we sample a query $q$ with
probability $\what{n}_q / n$, after which we uniformly sample one of the
${\what{n}_q\choose k}$ collections $\{i_1, \ldots, i_k\}$ of
$k$ indices
associated with query $q$, and we then perform the gradient
update \eqref{eqnfobos-update} using the gradient sample $\nabla
\surrloss(f^t, \structure(\preflabel_{i_1}, \ldots, \preflabel
_{i_k}))$. This
stochastic gradient scheme means that we can minimize the empirical
risk in a
number of iterations independent of both $n$ and $k$; the run-time
behavior of
the method scales independently of $n$ and depends on $k$ only so much as
computing an instantaneous gradient $\nabla\surrloss(f,
\structure(\preflabel_1, \ldots, \preflabel_k))$ increases with $k$.

\begin{figure}

\includegraphics{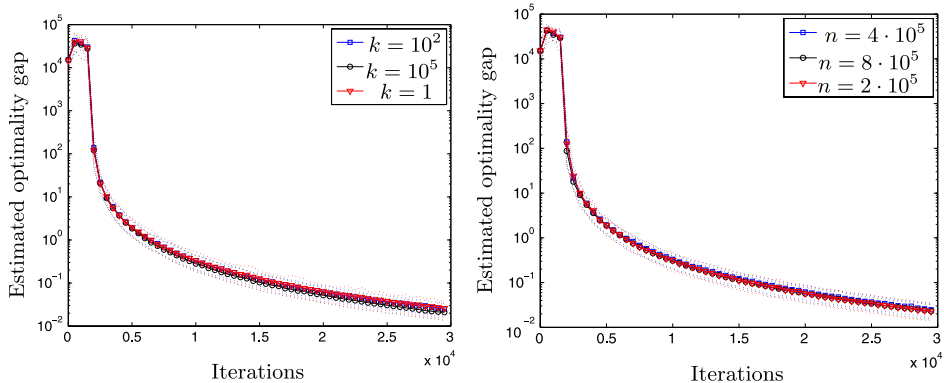}

\caption{Timing experiments for different values of $k$
and $n$ when applying the method \protect\eqref{eqnfobos-update}.
The horizontal
axes are the number of stochastic gradient iterations; the
vertical axes are the estimated optimality gap for the empirical surrogate
risk. Left: varying amount of aggregation $k$, fixed $n = 4 \cdot
10^5$. Right: varying total number of samples $n$, fixed $k =
10^2$.}\label{figtiming}
\end{figure}
%
%

In Figure~\ref{figtiming}, we show empirical evidence that the stochastic
method \eqref{eqnfobos-update} works as described. In particular, we minimize
the empirical $U$-statistic-based risk \eqref{eqnu-statistic-emprisk} with
the loss \eqref{eqnregression-loss} we employ in our experiments in
the next
section. In each plot in Figure~\ref{figtiming}, we give an estimated
optimality gap, $\emprisk_{\surrloss,n}(f^t) - \inf_{f \in\funclass}
\emprisk_{\surrloss,n}(f)$, as a function of $t$, the number of iterations.
As in the section to follow, $\funclass$ consists of linear functionals
parameterized by a vector $\theta\in\R^d$ with $d = 136$. To estimate
$\inf_{f \in\funclass}$, we perform\vspace*{1pt} 100,000 updates of the
procedure \eqref{eqnfobos-update}, then estimate $\inf_{f \in
\funclass}
\emprisk_{\surrloss,n}(f)$ using the output predictor $\what{f}$
evaluated on
an additional (independent) 50,000 samples (the number of terms in the true
objective is too large to evaluate). To estimate the risk
$\emprisk_{\surrloss,n}(f^t)$, we use a moving average of the
previous 100
sampled losses $\surrloss(f^\tau, \structure^{i_\tau})$ for $\tau
\in\{t -
99, \ldots, t\}$, which is an unbiased estimate of an upper bound on the
empirical risk $\emprisk_{\surrloss, n}(f^t)$ (see,
e.g., \cite{CesaBianchiCoGe02}). We perform the experiment 20 times
and plot
averages as well as 90\% confidence intervals. As predicted by our theoretical
results, the number of iterations to attain a particular accuracy is
essentially independent of $n$ and $k$; all the plots lie on one another.

\subsection{Experimental evaluation}
\label{secmsr-evaluation}

To perform our experimental evaluation, we use a subset of the Microsoft
Learning to Rank Web10K dataset \cite{QinLiDiXuLi11}, which consists of
10,000 web searches (queries) issued to the Microsoft Bing search
engine, a
set of approximately 100 potential results for each query, and a relevance
score $r \in\R$ associated with each query/result pair. A
query/result pair
is represented by a $d = 136$-dimensional feature vector of standard
document-retrieval features.

To understand the benefits of aggregation and consistency in the
presence of
partial preference data, we generate pairwise data from the observed
query/result pairs, so that we know the true asymptotic generating
distribution. We adopt a loss $L$ from the NDCG-family \eqref
{eqnndcg} and
compare three surrogate losses: a Fisher-consistent regression
surrogate based
on aggregation, an inconsistent but commonly used pairwise logistic
loss \cite{DekelMaSi03}, and a Fisher-consistent loss that requires
access to
complete preference data \cite{RavikumarTeYa11}. Recalling the NDCG
score \eqref{eqnndcg} of a prediction vector $\alpha\in\R
^\numresults$ for
scores $\structure\in\R^\numresults$ (where $\pi_\alpha$ is the
permutation
induced by~$\alpha$), we have the loss
\[
L(\alpha, \structure) = 1 - \frac{1}{Z(\structure)} \sum_{j=1}^\numresults
\frac{G(\structure_j)}{F(\pi_\alpha(j))},
\]
where $Z(\structure)$ is the normalizing value for the NDCG
score, and $F(\cdot)$ and $G(\cdot)$ are increasing functions.

Given a set of queries $q$ and relevance scores $r_i \in\R$,
we generate $n$ pairwise preference observations according to a
Bradley--Terry--Luce (BTL) model \cite{BradleyTe52}.
That is, for each observation, we choose a query $q$ uniformly at random
and then select a uniformly random pair $(i,j)$ of results to compare.
The pair is ordered as $i \succ j$ (item $i$ is preferred to $j$)
with probability $p_{ij}$, and $j \succ i$ with probability $1 - p_{ij} =
p_{ji}$, where
%
%
\begin{equation}
\label{eqnlogistic-pairwise-sampling} p_{ij} = \frac{\exp(r_i - r_j)}{1 + \exp(r_i - r_j)}
\end{equation}
for $r_i$ and $r_j$ the respective relevances of results $i$ and $j$ under
query $q$.

We define our structure functions $\structure_k$ as score vectors in
$\R^\numresults$, where given a set of $k$ preference pairs, the
score for
item $i$ is
\[
\structure_k(i) = \frac{1}{\numresults- 1} \sum
_{j \neq i} \log\frac{\what{\P}(j \prec i)}{\what{\P}(j \succ i)},
\]
the average empirical log-odds of result $i$ being preferred to any other
result. Under the BTL model \eqref{eqnlogistic-pairwise-sampling}, as
$k\to\infty$ the structural score converges\vadjust{\goodbreak} for\break each $i \in
[\numresults]$ to
%
%
\begin{equation}
\structure(i) = \frac{1}{\numresults- 1} \sum_{j \neq i}
\bigl[\log\bigl(1 + \exp(r_i - r_j)\bigr) - \log\bigl(1 +
\exp(r_j - r_i)\bigr) \bigr]. \label{eqnexperiment-limiting-score}
\end{equation}
In our setting, we may thus evaluate the asymptotic NDCG risk
of a scoring function $f$ by computing the asymptotic
scores \eqref{eqnexperiment-limiting-score}. In addition,
Corollary~\ref{corollaryndcg-consistency} shows that if all
minimizers of a loss obey the ordering of the values
\[
\int_\structurespace \frac{G(\structure(j))}{Z(\structure)} \,d\limitlaw(\structure),\qquad j \in
\{1, \ldots, \numresults\}
\]
then the loss is Fisher-consistent. A well-known
example \cite{CossockZh08,RavikumarTeYa11} of such a loss is the
least-squares loss, where the regression labels are $G(\structure(j)) /
Z(\structure)$:
%
%
\begin{equation}
\surrloss(\alpha, \structure) = \frac{1}{2 \numresults} \sum
_{j = 1}^\numresults \biggl(\alpha_j -
\frac{G(\structure(j))}{Z(\structure)} \biggr)^2. \label{eqnregression-loss}
\end{equation}
We compare the least-squares aggregation loss with a pairwise
logistic loss natural for the pairwise data generated according to the
BTL model \eqref{eqnlogistic-pairwise-sampling}. Specifically,
given a data pair with $i \succ j$, the logistic surrogate loss is
%
%
\begin{equation}
\label{eqnlogistic-loss} \surrloss(\alpha, i \succ j) = \log \bigl(1 + \exp(
\alpha_j - \alpha_i) \bigr),
\end{equation}
which is equivalent or similar to previous losses used for pairwise
data in
the ranking literature \cite{Joachims02,DekelMaSi03}. For
completeness, we
also compare with a Fisher-consistent surrogate that requires access to
complete preference information in the form of the asymptotic structure
scores \eqref{eqnexperiment-limiting-score}. Following
Ravikumar et~al. \cite{RavikumarTeYa11}, we obtain such a surrogate by granting the regression
loss \eqref{eqnregression-loss} direct access to the asymptotic structure
scores. Note that such a construction would be infeasible in any true
pairwise data setting.

Having described our sampling procedure, aggregation strategy, and loss
functions, we now describe our model. We let $x^q_i$ denote the feature vector
for the $i$th result from query $q$, and we model the scoring function $f(q)_i
=  \langle \theta, x^q_i \rangle $ for a vector $\theta
\in\R^d$. For the regression
loss \eqref{eqnregression-loss}, we minimize the $U$-statistic-based
empirical risk \eqref{eqnu-statistic-emprisk} over a variety of
orders $k$,
while for the pairwise logistic loss~\eqref{eqnlogistic-loss}, we minimize
the empirical risk over all pairs sampled according to the BTL
model~\eqref{eqnlogistic-pairwise-sampling}. We regularize our
estimates by
adding $\regularizer(\theta) = (\lambda/2)\llVert {\theta}\rrVert
_2^2$ to the objective
minimized, and we use the specialized stochastic
method \eqref{eqnfobos-update} to minimize the empirical risk.

\begin{figure}

\includegraphics{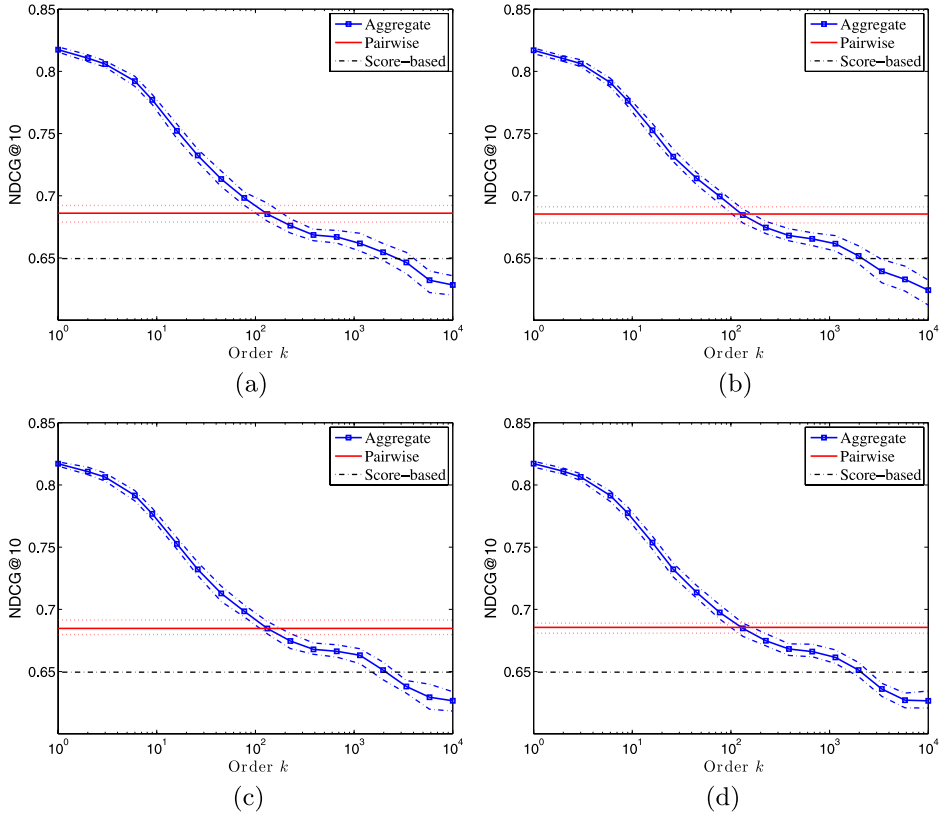}
\vspace*{-3pt}
\caption{NDCG risk and 95\% confidence intervals for
$\theta$ estimated using the logistic pairwise
loss~\protect\eqref{eqnlogistic-loss} and the $U$-statistic
empirical risk with
$\surrloss$ chosen to be regression
loss \protect\eqref{eqnregression-loss}. The horizontal axis of each
plot is
the order $k$ of the aggregation in the
$U$-statistic \protect\eqref{eqnu-statistic-emprisk}, the vertical
axis is the
NDCG risk, and each plot corresponds to a different number $n$ of
samples. \textup{(a)} $n = 2 \cdot10^5$; \textup{(b)} $n = 4 \cdot10^5$; \textup{(c)} $n = 8 \cdot
10^5$; \textup{(d)} $n = 1.6 \cdot10^6$.}\label{figk-orders}\vspace*{-5pt}
\end{figure}
%
%
%
%
%
%
%

\begin{figure}

\includegraphics{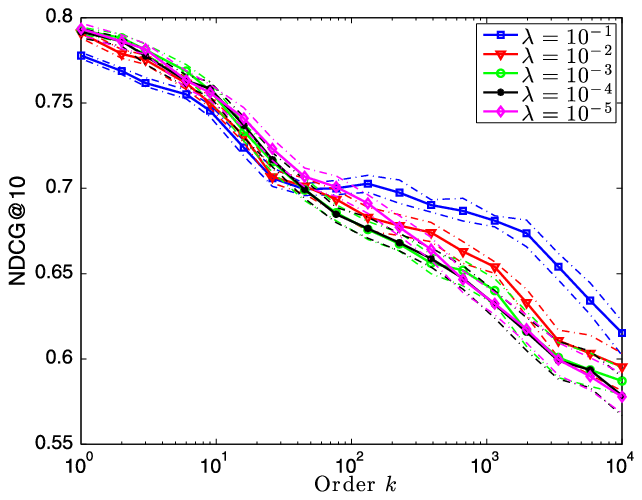}

\caption{NDCG risk and 95\% confidence intervals
for $\theta$ estimated using the $U$-statistic empirical
risk~\protect\eqref{eqnu-statistic-emprisk} with $\surrloss$ chosen
as the
regression loss \protect\eqref{eqnregression-loss} under various
choices of the
regularization parameter, $\lambda$.}
\label{figlambda-size}
\end{figure}

%

Our goals in the experiments are to understand the behavior of the empirical
risk minimizer as the order $k$ of the aggregating statistic is varied
and to
evaluate the extent to which aggregation improves the estimated
scoring function. A secondary concern is to verify that the
method is insensitive to the amount $\lambda$ of regularization performed
on $\theta$. We run each experiment 50 times and report
confidence intervals based on those 50 experiments.

Let $\theta^\mathrm{reg}_{n,k}$ denote the estimate of $\theta$ obtained from
minimizing the empirical risk~\eqref{eqnu-statistic-emprisk} with the
regression loss \eqref{eqnregression-loss} on $n$ samples with aggregation
order $k$, let $\theta^\mathrm{log}_n$ denote the estimate of $\theta$ obtained
from minimizing the empirical pairwise logistic
loss \eqref{eqnlogistic-loss}, and let $\theta^\mathrm{full}$ denote the
estimate of $\theta$ obtained from minimizing the empirical risk with surrogate
loss \eqref{eqnregression-loss} using the asymptotic structure
scores~\eqref{eqnexperiment-limiting-score} directly. Then each plot of
Figure~\ref{figk-orders} displays the risk $\risk(\theta^\mathrm{reg}_{n,k})$ as
a function of the aggregation order $k$, using $\risk(\theta^\mathrm{log}_n)$ and
$\risk(\theta^\mathrm{full})$ as references. The four plots in the figure
correspond to different numbers $n$ of data pairs.

Broadly, the four plots in Figure~\ref{figk-orders} match our theoretical
results. Consistently across the plots, we see that for small $k$, it appears
there is not sufficient aggregation in the regression-loss-based empirical
risk, and for such small $k$ the\vadjust{\goodbreak} pairwise logistic loss is better.
However, as
the order of aggregation $k$ grows, the risk performance of $\theta
_{n,k}^\mathrm{reg}$ improves. In addition, with larger sample sizes $n$, the difference
between the risk of $\theta^\mathrm{log}_n$ and $\theta^\mathrm{reg}_{n,k}$ becomes
more pronounced. The second salient feature of the plots is a moderate
flattening of the risk $\risk(\theta^\mathrm{reg}_{n,k})$ and widening
of the
confidence interval for large values of $k$. This seems consistent with the
estimation error guarantees in the
theoretical results in Lemmas~7
and~10 in the appendices,
where the order $k$ being
large has an evidently detrimental effect. Interestingly, however, large
values of $k$ still yield significant improvements over $\risk(\theta
^\mathrm{log}_n)$. For very large $k$, the improved performance of $\theta^\mathrm{reg}_{n,k}$ over $\theta^\mathrm{full}$ is a consequence of sampling artifacts
and the fact that we use a finite dimensional representation. [By
using sufficiently many dimensions $d$, the estimator $\theta^\mathrm{full}$
attains zero risk by matching the asymptotic
scores \eqref{eqnexperiment-limiting-score} directly.]

Figure~\ref{figlambda-size} displays the risk $\risk(\theta^\mathrm{reg}_{n,k})$
for $n = 800\mbox{,}000$ pairs, $k = 100$, and multiple values of the
regularization multiplier $\lambda$ on $\llVert {\theta}\rrVert
_2^2$. The
results, which
are consistent across many choices of $n$, suggest that minimization of the
aggregated empirical risk \eqref{eqnu-statistic-emprisk} is robust to the
choice of regularization multiplier.


\section{Conclusions}
\label{secconclusions}

In this paper, we demonstrated both the difficulty and the
feasibility of designing consistent, practicable procedures for ranking.
By giving necessary and sufficient conditions for the Fisher
consistency of
ranking algorithms, we proved that many natural ranking procedures
based on
surrogate losses are inconsistent, even in low-noise settings. To address
this inconsistency while accommodating the incomplete nature of typical
ranking data, we proposed a new family of surrogate losses, based on
$U$-statistics, that aggregate disparate partial preferences. We showed how
our losses can fruitfully leverage\vadjust{\goodbreak} any well behaved rank aggregation procedure
and demonstrated their empirical benefits over more standard surrogates
in a
series of ranking experiments.

Our work thus takes a step toward bringing the consistency literature for
ranking in line with that for classification, and we anticipate several
directions of further development. First, it would be interesting to
formulate low-noise conditions under which faster rates of convergence are
possible for ranking risk minimization (see, e.g.,
the work of \cite{ClemenconLuVa08}, which focuses on the
minimization of a single pairwise loss). Additionally, it
may be interesting to study structure functions $\structure$ that yield
nonpoint distributions $\limitlaw$ as the number of arguments $k$
grows to
infinity. For example, would scaling the Thurstone--Mosteller least-squares
solutions \eqref{eqnskew-symmetric-score} by $\sqrt{k}$---to achieve
asymptotic normality---induce greater robustness in the empirical
minimizer of
the $U$-statistic risk \eqref{eqnu-statistic-emprisk}? Finally, exploring
tractable formulations of other supervised learning
problems in which label data is naturally incomplete could be fruitful.

%
%


\section*{Acknowledgments}
We thank the anonymous reviewers and the Associate Editor for their
helpful comments and valuable feedback.

\begin{supplement}[id=suppA]
\stitle{Proofs of results}
\slink[doi]{10.1214/13-AOS1142SUPP} 
\sdatatype{.pdf}
\sfilename{aos1142\_supp.pdf}
\sdescription{The supplementary material contains proofs of our results.}
\end{supplement}

%

\printaddresses

\end{document}